\documentclass[a4paper,11pt]{article}

\usepackage[utf8]{inputenc} 
\usepackage{amsfonts}
\usepackage{amssymb}
\usepackage{amsthm}
\usepackage{amsmath}
\usepackage{esint} % For the average integral
\usepackage{stackengine}
\usepackage{a4wide}
\usepackage{mathrsfs}
\usepackage{epsfig}
\usepackage{palatino}
\usepackage{tikz}
\usepackage{fp,ifthen}
\usepackage{nicefrac}
\usetikzlibrary{decorations.pathreplacing}
\usepackage{float}
\usepackage{enumerate} %Enumerate package must be always written before enumitem
\usepackage{enumitem} %Customize enumerate and itemize
\usepackage{comment}

\usepackage[a4paper,twoside,top=2.5cm, bottom=2cm, left=2.3cm, right=2.3cm]{geometry} 
 \usepackage{lipsum}

\usepackage[colorinlistoftodos,textwidth=2.3cm]{todonotes} %Package for comments, WARNING: it does not compile if you write Italian stresses like è,ì,à,ò,.. inside the command \todo{...}

\newcommand{\pdfgraphics}{\ifpdf\DeclareGraphicsExtensions{.pdf,.jpg}\else\fi}
\usepackage{graphicx}
\usepackage{dsfont}

\usepackage{color}
\definecolor{hanblue}{rgb}{0.27, 0.42, 0.81}
\definecolor{red}{rgb}{1.0, 0.0, 0.0}
\usepackage[colorlinks, citecolor=blue,linkcolor=blue, urlcolor = blue]{hyperref}
\usepackage[english,capitalize]{cleveref}

\usepackage{mathtools}
\theoremstyle{plain}

\newtheorem{teo}{Theorem}[section]
\newtheorem{lemma}[teo]{Lemma}
\newtheorem{prop}[teo]{Proposition}
\newtheorem{cor}[teo]{Corollary}
\newtheorem{ackn}{Acknowledgements\!}

\theoremstyle{definition}
\newtheorem{defn}[teo]{Definition}

\theoremstyle{remark}
\newtheorem{remark}[teo]{Remark}

\numberwithin{equation}{section}

\newcommand{\de}{\ensuremath{\,\mathrm d}} % Il de degli integrali (contiene lo spazio da mettere tra integranda e misura)
\newcommand \eps{\ensuremath{\varepsilon}} 
\renewcommand{\epsilon}{\varepsilon}
\newcommand{\N}{\ensuremath{\mathbb N}}%Natural numbers
\newcommand{\R}{\ensuremath{\mathbb R}}%%Real numbers
\newcommand{\X}{\mathsf{X}} % Metrix space X
\newcommand{\Y}{\mathsf{Y}} % Metric space Y
\newcommand{\Z}{\mathsf{Z}}%Metric space Z
\newcommand{\dist}{{\rm d}} % Distance function.
\newcommand{\lip}{{\rm lip}} %Local Lipschitz constant
\newcommand{\m}{\ensuremath{\mathfrak m}}

\DeclareMathOperator\Lip{Lip}

\DeclareMathOperator*{\argmin}{arg\,min}

\makeatletter
\renewcommand*\env@matrix[1][*\c@MaxMatrixCols c]{%
  \hskip -\arraycolsep
  \let\@ifnextchar\new@ifnextchar
  \array{#1}}
\makeatother

\begin{document}

\pdfgraphics % Use this command right after \begin{document}

\title{Comments on the regularity of harmonic maps between singular spaces}

\author{Luca Gennaioli \footnote{\href{luca.gennaioli@sissa.it}{luca.gennaioli@sissa.it}} \and Nicola Gigli \footnote{\href{nicola.gigli@sissa.it}{nicola.gigli@sissa.it}} \and Hui-Chun Zhang
\footnote{\href{zhanghc3@mail.sysu.edu.cn}{zhanghc3@mail.sysu.edu.cn}} \and Xi-Ping Zhu
\footnote{\href{stszxp@mail.sysu.edu.cn}{zstszxp@mail.sysu.edu.cn}}}

\date{\today}

\maketitle

\vspace{-0.5cm}
\begin{abstract}
\noindent 
In this work we are going to establish H\"older continuity of harmonic maps from an open set $\Omega$ in an ${\rm RCD}(K,N)$ space valued into a ${\rm CAT}(\kappa)$ space, with the constraint that the image of $\Omega$ via the map is contained in a sufficiently small ball in the target. Building on top of this regularity and assuming a local Lipschitz regularity of the map, we establish a weak version of the Bochner-Eells-Sampson inequality in such a non-smooth setting. Finally we study the boundary regularity of such maps.
\end{abstract}

\tableofcontents

\section{Introduction}
In the last 50 years the study of harmonic maps has been blooming and gained a lot of interest from the mathematical community. One of the main questions is the one of existence of such mappings and parallel to that there is the issue of their regularity.

When $u:\Omega\subseteq M^n\to N^k$ is an harmonic map between Riemannian manifolds $(M^n,g_M)$ and $(N^k,g_N)$ the picture nowadays is quite clear: the existence of such mappings has been established via the study of parabolic problems by Hamilton (see \cite{Hamilton75} for a discussion on the topic) and then by looking at the problem in a variational way. The latter approach can be tailored to the non-smooth setting as well, indeed in the recent \cite{Sak23} the author has been able to prove the existence of harmonic maps between an ${\rm RCD}$ space and a ${\rm CAT}(\kappa)$ space if the image of $u$ is contained in a sufficiently small ball. 

Back to the case of an harmonic map between smooth Riemannian manifolds, the Bochner-Eells-Sampson formula states that 
\begin{equation*}
    \Delta\bigg(\frac{|\de u|_{\sf HS}^2}{2}\bigg)=|\nabla\de u|_{\sf HS}^2+{\rm Ric}_{g_M}(\nabla u,\nabla u)-\sum_{i,j\leq n}\langle u_*\mathcal{R}^N(e_i,e_j)e_i,e_j\rangle,
\end{equation*}
where ${\rm Ric}_{g_M}$ is the Ricci tensor of the source space, $u_*\mathcal{R}^N$ is the pullback of the curvature tensor of the target space via the map $u$ and ${\rm e_\alpha}_{\alpha=1}^n$ is an orthonormal frame for the tangent bundle $TM$. If we assume that $Ric_{g_M}\geq -K$ (lower bound on the Ricci tensor) and ${\rm R_N}\leq \kappa$ (upper bound on the sectional curvatures), the previous identity can be turned into the following inequality
\begin{equation}
    \label{Smooth BES inequality}
      \Delta\bigg(\frac{|\de u|_{\sf HS}^2}{2}\bigg)\geq |\nabla\de u|_{\sf HS}^2+ K|\de u|_{\sf HS}^2-\kappa |\de u|_{\sf HS}^4.
\end{equation}
From this inequality, at least if $\kappa=0$ it is possible to quickly deduce that harmonic maps are locally Lipschitz, as in this case we have
\begin{equation}
\label{eq:kappazero}\Delta\bigg(\frac{|\de u|_{\sf HS}^2}{2}\bigg)\geq K|\de u|_{\sf HS}^2
\end{equation}
and thus a De Giorgi-Nash-Moser argument shows that the function $f:=|\de u|_{\sf HS}^2$ is locally bounded. The case $\kappa>0$, say $\kappa=1$, is more delicate and is known to require the additional assumption that the range of $u$ is contained in a ball $B_r(p)\subset N^k$ of radius $r<\tfrac\pi2$ (otherwise there are known counterexamples to regularity \cite{Riv95}). On top of this,  the term $|\de u|^4$ is a priori not in $L^1$, making  it hard to extract information from \eqref{Smooth BES inequality}. To overcome these difficulties,  Serbinowski argued as follows: the function $f(x):={\sf d}_N(u(x),p)$ satisfies $-\Delta\cos(f)\geq |\de u|^2_{\sf HS}\cos(f)$ (as a consequence of the fact that $u$ is harmonic and of the curvature assumption on $N$), and quite trivially we have $|\de|\de u|_{\sf HS}|^2\leq |\nabla\de u|_{\sf HS}^2$. These consideration and little algebraic manipulation show that \eqref{Smooth BES inequality} implies
\begin{equation}
    \label{Serb}
    \tfrac{|\de u|_{\sf HS}}{\cos(f)}{\rm div}\Big(\cos^2(f)\nabla\Big(\tfrac{|\de u|_{\sf HS}}{\cos(f)}\Big)\Big)\geq K|\de u|^2_{\sf HS},
\end{equation}
and since $\cos(f)$ is far from zero, a Moser iteration argument can be called into play to prove that $\tfrac{|\de u|_{\sf HS}}{\cos(f)}$, and thus $|\de u|_{\sf HS}$, is locally bounded, as desired.

This type of reasoning allows to conjecture that, in the non-smooth setting, one should impose a lower bound on the Ricci curvature of the source and an upper bound on the sectional curvature on the target to get that an harmonic maps is (locally) Lipschitz.

Many contributions in this direction have appeared in the recent years: for an account of the story we refer to the extensive introductions in \cite{ZZ18}, \cite{ZZZ19}, \cite{MS22} and \cite{Gig23}. Here  we just recall that in  \cite{ZZ18} the authors proved the Lipschitz regularity of harmonic maps between Alexandrov spaces and a weak Bochner-Eells-Sampson inequality. Building on this, in the more recent \cite{Gig23} and \cite{MS22} the authors where able to establish such regularity when the source space is an ${\rm RCD}(K,N)$ space, namely a space with a synthetic notion of Ricci curvature bounded below by $K$ and dimension bounded above by $N$, and the target is a ${\rm CAT}(0)$ space, namely a space with a synthetic notion of sectional curvature bounded above by $0$.

Very roughly said, the basic argument to get a sort of \eqref{eq:kappazero} and local Lipschitz regularity of harmonic maps is to build two families $(g_t),(h_t)$ of functions  (via a kind of Hopf-Lax formula for metric-valued maps) converging to $|\de u|^2$ in $L^1$ as $t\downarrow0$ satisfying
\begin{equation}
\tfrac12\Delta g_t\geq K h_t\qquad\forall t>0.
\end{equation}
Quite clearly, from this it is possible to pass to the limit and obtain that
\begin{equation}
\label{eq:zzz}
\Delta\bigg(\frac{|\de u|^2}{2}\bigg)\geq K|\de u|^2.
\end{equation}
Notice that in this the quantity $|\de u|$ is the operator norm of $\de u$, not its Hilbert-Schmidt norm as in \eqref{eq:kappazero}, thus  \eqref{eq:zzz} is not the same as \eqref{eq:kappazero}, but the effect is the same: a Moser iteration argument shows that $|\de u|$ must be locally bounded and thus that $u$ is locally Lipschitz.

When dealing with the case $\kappa>0$ this strategy encounters a problem, as the approximation procedure does not work well in conjunction with Serbinowski's technique: shortly said, at the approximated level the right hand side of \eqref{Smooth BES inequality} still contains a term that does not go to 0 in $L^1$ as $t\downarrow0$.

Because of these difficulties, we do not achieve Lipschitz regularity of harmonic maps in the mote general setting, our main results are rather:
\begin{itemize}
\item[1)] the proof of H\"older continuity, see Theorem \ref{Holder continuity}. Here we follow the strategy in \cite{Jost1997}.
\item[2)] the higher integrability of the energy density, see Theorem \ref{higher-integrable}, by using a Caccioppoli inequality and the Gehring lemma in \cite{Maasalo07,AndersonTH17}.
\item[3)] Under the a priori assumption that the harmonic function is Lipschitz, possibly with a sub-optimal control on the Lipschitz constant, we prove a version of inequality \eqref{Smooth BES inequality}, see Theorem \ref{thm:varbes}. To achieve this we suitable combine ideas from \cite{ZZZ19}, \cite{MS22} and \cite{Gig23}. Once we have this,  following the arguments in \cite{ZZZ19} one can obtain a sharp estimate on the Lipschitz constant and, as a consequence, a Liouville-type of result,  Theorem \ref{thm:sharplip} and Corollary \ref{Yauineq} for the precise statements.
\item[4)] the boundary regularity, see Theorem \ref{thm-boundary-reg}.  
\end{itemize}
\begin{ackn}
    N.G. acknowledges the support of the European Union - NextGenerationEU, in the framework of the PRIN Project 'Contemporary perspectives on geometry and gravity' (code 2022JJ8KER – CUP G53D23001810006). The views and opinions expressed are solely those of the authors and do not necessarily reflect those of the European Union, nor can the European Union be held responsible for them. H.C.Z. is partially supported by NSFC 12025109, and X.P.Z. is partially supported by NSFC 12271530. 
\end{ackn}

\section{Preliminaries}
\subsection{The source: ${\rm RCD}(K,N)$ spaces}
We say that $(\X,\dist_\X,\m)$ is a metric measure space if $(\X,\dist_\X)$ is a complete and separable metric space and $\m$ is a Radon measure which is finite on balls. For a function $f:\X\to\R$ we set 

\begin{equation*}
    \lip f(x):=\begin{cases}
    \limsup_{y\to x}\frac{|f(y)-f(x)|}{\dist_\X(x,y)}\qquad&{\rm if\;x\;is\;not\;isolated} \\ 
    0\qquad&{\rm if\;x\;is\;isolated}
    \end{cases}
\end{equation*}
and we call it \emph{local Lipschitz constant} of $f$, while with $\Lip f$ we denote the classical Lipschitz constant of $f$.

In order to develop Sobolev calculus on metric measure space, following \cite{Cheeger00}, we introduce the Cheeger energy ${\rm Ch}:L^2(\m_\X)\to[0,\infty]$ as 
\begin{equation*}
    {\rm Ch}(f):=\inf\bigg\{\liminf_{k\to\infty}\frac{1}{2}\int_{\X}\lip^2(f_k)\de\m_\X:\; (f_k)_k\in\Lip_{\rm bs}(\X),\;f_k\to f\;{\rm in}\;L^2(\m_\X)\bigg\}.
\end{equation*}
It can be proved that if ${\rm Ch}(f)<\infty$ there exists a function, which we call $|\nabla f|$, such that $|\nabla f|\in L^2(\m_\X)$ and 
\begin{equation*}
    {\rm Ch}(f)=\frac{1}{2}\int_\X|\nabla f|^2\de\m_\X.
\end{equation*}
If that is the case we say that $f\in W^{1,2}(\X)$. The latter set is actually a vector space with its natural operation and, if endowed with the norm $\|f\|_{W^{1,2}}=\|f\|_{L^2}+2{\rm Ch}(f)$, it is also a Banach space. In order to introduce a well-behaved notion of Laplacian of a function we shall now speak about infinitesimal Hilbertianity.
We say that a metric measure space is \emph{infinitesimally Hilbertian}, following \cite{Gigli12}, if ${\rm Ch}$ is a quadratic form. In this case via polarization it is possible to give a meaning to the object 
\begin{equation*}
    \int_\X\langle\nabla\varphi,\nabla f\rangle\de\m_\X
\end{equation*}
by setting 
\begin{equation*}
   \int_\X\langle\nabla\varphi,\nabla f\rangle\de\m_\X:= {\rm Ch}(f+\varphi)- {\rm Ch}(f)- {\rm Ch}(\varphi)
\end{equation*}
\begin{defn}[$L^2$ Laplacian]
    We say that $f:\X\to\R$ in $W^{1,2}(\X)$ is such that $f\in D(\Delta)\subset L^2(\m_\X)$ if there exists $g\in L^2(\m_\X)$ such that 
    \begin{equation*}
         -\int_\X\langle\nabla\varphi,\nabla f\rangle\de\m_\X=\int_\X g\varphi\de\m_\X
    \end{equation*}
    for all $\varphi\in W^{1,2}(\X)$. We shall set $\Delta f:=g$.
\end{defn}
\begin{defn}[Measure-valued Laplacian]
    We say that $f:\X\to\R$ in $W^{1,2}_{\rm loc}(\X)$ has measure-valued Laplacian in $\Omega$ if there exists a Radon measure $\mu\in\mathcal{M}(\Omega)$ such that 
    \begin{equation*}
        -\int_\X\langle\nabla\varphi,\nabla f\rangle\de\m_\X=\int_\X\varphi\de\mu
    \end{equation*}
    for all $\varphi\in \Lip_c(\Omega)$, the latter being the space of Lipschitz functions with compact support inside $\Omega$.
\end{defn}
\begin{remark}
    With a little bit of abuse of notation we shall call $\Delta f=\mu$ the measure-valued Laplacian as well. We will do this since if $\mu\ll \m_\X$ with density in $L^2_{\rm loc}$, then $\mu=\Delta f\m_\X$. Notice also that we are using the term \emph{Radon measures} to denote what are more properly called Radon functionals (see \cite{CavMon20}).
\end{remark}

We are now ready to introduce the class of spaces which we will use as source space for the definition of our harmonic map $u$. 
We can introduce ${\rm RCD}(K,N)$ spaces building on the tools we have just presented. Following an Eulerian approach it is possible to characterize them via the Bochner inequality (see \cite{Gigli-Kuwada-Ohta10}, \cite{AmbrosioGigliSavare11-2}, \cite{AmbrosioGigliSavare12}, \cite{Erbar-Kuwada-Sturm13}, \cite{AmbrosioMondinoSavare13}, \cite{CavMil16}). For a more detailed discussion on such notions and for the interplay with optimal transport we refer to the recent \cite{gigli2023giorgi} and \cite{AmbICM}.
\begin{defn}[${\rm RCD}(K,N)$ space]
We say that a metric measure space $(\X,\dist_\X,\m_\X)$ is an ${\rm RCD}(K,N)$ space if the following conditions are met:
\begin{enumerate}
    \item There exists $c_1,c_2\geq0$ such that for some $x\in\X$ we have
    \begin{equation*}
        \m(B_r(x))\leq C_1 e^{c_2 r^2}.
    \end{equation*}
    \item $W^{1,2}(\X)$ is a Hilbert space. \\ 
    \item If $f\in W^{1,2}(\X)$ is such that $|\de f|\leq 1$ $\m$-a.e., then $f$ has a $1$-Lipschitz representative. \\
    \item For every $f\in D(\Delta)$ with $\Delta f\in W^{1,2}(\X)$ and $g\in L^\infty(\m)\cap D(\Delta)$ the following \emph{Bochner inequality} holds 
    \begin{equation*}
        \int_\X\frac{|\de f|^2}{2}\Delta g\de\m\geq \int_\X g\bigg(K|\de f|^2+\frac{(\Delta f)^2}{N} + \langle\nabla f,\nabla\Delta f\rangle\bigg)\de\m.
    \end{equation*}
\end{enumerate}
\end{defn}
The final object we shall introduce is the heat semigroup $h_t:L^2(\m_\X)\to L^2(\m_\X)$: it can be introduced as the gradient flow of the Cheeger energy. Therefore we shall call $(h_tf)_{t\geq 0}$ such a gradient flow starting from $f\in L^2(\m_\X)$. For an account of its properties the reader can consult \cite{GP19}. If the space $\X$ is an ${\rm RCD}(K,N)$ space then it is possible to consider the ${\rm EVI}_K$ gradient flow of the entropy functional on the space of probability measures. If we denote with $h_t\delta_x$ the gradient flow of the entropy starting from a Dirac mass centered at $x$ we have $h_t\delta_x\ll \m_\X$ and we shall call $p_t(x,y):=\frac{\de h_t\delta_x}{\de\m_\X}(y)$. It can be proved that $h_tf:=\int_\X p_t(x,\cdot)f(x)\de\m_\X$ and that $p_t$ is H\"older continuous and satisfies the following Gaussian estimates 
\begin{equation}
    \label{Gaussian estimates}
  \frac{c}{\m_\X(B_{\sqrt{t}}(x))}e^{-\dist_\X^2(x,y)/3t-C_1t}\leq p_t(x,y)\leq \frac{C}{\m_\X(B_{\sqrt{t}}(x))}Ce^{-\dist_\X^2(x,y)/5t+C_2t},
\end{equation}
for all $x,y\in\X$, $t>0$ and for some $c,C,C_1,C_2>0$.
There is also a gradient bound thanks to the Li-Yau inequality but for the sake of exposition we shall limit ourselves to this presentation: the interested reader can consult \cite{JLZ16}, \cite{KTS94-1}, \cite{KTS94-2} and \cite{KTS94-3} for more information on Gaussian estimates.

Since we are interested in giving a meaning to "$\Delta f\geq \eta$" we shall rigorously introduce such a notion:
\begin{defn}[Weak Laplacian bound]
\label{Weak Laplacian bound}
    Let $(\X,\dist_\X,\m)$ be a metric measure space and $\Omega\subset\X$ an open and bounded set. Let $\eta:\Omega\to\R$ be continuous and bounded.  We say that a function $f\in W^{1,2}_{\rm loc}(\Omega)$ is such that $\Delta f\leq \eta$ in the weak sense if for all $\varphi\in\Lip_c^+(\Omega)$ (being $\Lip_c^+(\Omega)$ the subset of $\Lip_c(\Omega)$ made of nonnegative functions) we have 
    \begin{equation*}
        -\int_{\X}\nabla f\cdot\nabla\varphi\de\m_\X\leq\int_\X\varphi\eta\de\m_\X.
    \end{equation*}
\end{defn}
To introduce another (weaker) notion of Laplacian bounds we need to introduce the following space
\begin{equation*}
    {\rm Test}_c^\infty(\X):=\bigg\{\varphi\in D(\Delta)\cap L^\infty: |\nabla\varphi|\in L^\infty , \Delta\varphi\in L^\infty\cap W^{1,2}\bigg\}.
\end{equation*}
We write ${\rm Test}_c^\infty(\Omega)$ if ${\rm supp}\varphi\subset\subset\Omega$.
\begin{comment}
\begin{defn}[Very weak Laplacian bound]
\label{Very weak Laplacian bound}
    Let $(\X,\dist_\X,\m_\X)$ be a metric measure space and $\Omega\subset\X$ an open and bounded set. Let $\eta:\Omega\to\R$ be continuous and bounded.  We say that a function $f\in L^1_{\rm loc}(\Omega)$ is such that $\Delta f\leq \eta$ in the very weak sense if for all nonnegative $\varphi\in{\rm Test}_c^\infty(\Omega)$ we have 
    \begin{equation*}
        \int_{\X}f\Delta\varphi\de\m_\X\leq\int_\X\varphi\eta\de\m_\X.
    \end{equation*}
\end{defn}
\begin{remark}
\label{weak and very weak equivalence}
    As it was pointed out in \cite{Gigli_2023}, if $\X$ is an ${\rm RCD}(K,N)$ space and $f\in W^{1,2}_{\rm loc}(\Omega)$ has a very weak Laplacian bound, then such a bound holds also in the weak sense due to the density of ${\rm Test}_c^\infty(\Omega)$ into $\Lip_c(\Omega)$.
\end{remark}
\end{comment}
\begin{defn}[Heat flow Laplacian bound]
Let $(\X,\dist_\X,\m_\X)$ be an infinitesimally Hilbertian metric measure space and $\Omega\subset\X$ be an open and bounded set. Let $f:\Omega\to\R$ be a bounded and lower semicontinuous function and let $\eta\in C_b(\Omega)$. We say that $\Delta f\leq\eta$ in the heat flow sense if 
\begin{equation*}
    \limsup_{t\to 0}\frac{h_t\Tilde{f}(x)-\Tilde{f}(x)}{t}\leq\eta(x)
\end{equation*}
    for all $x\in\Omega$, where $\Tilde{f}:\X\to\R$ is the global extension of $f$ which is set to zero outiside of $\Omega$.
\end{defn}
Finally we recall the classical Laplacian comparison for the distance function from a point, which in this non-smooth setting has been obtained in \cite[Corollary 5.15]{Gigli12}.
\begin{teo}[Laplacian comparison]
\label{Laplacian comparison}
Let $(\X,\dist_\X,\m_\X)$ be an ${\rm RCD}(K,N)$ space for some $K\in\R$, $N\in\N$ and fix $x_0\in\X$. Then the map $x\to\dist^2_\X(x_0,x)=\dist^2_{\X,x_0}(x)$ has measure-valued Laplacian and 
\begin{equation*}
    \Delta\frac{\dist^2_{\X,x_0}}{2}\leq C(N,K,\dist_{\X,x_0}(\cdot))\m_\X
\end{equation*}
    in the weak sense. Moreover the same holds for the map $x\to\dist_{\X,x_0}(x)$, on $\X\setminus\{x_0\}$, namely
    \begin{equation*}
        \Delta\dist_{\X,x_0 |\X\setminus\{x_0\}}\leq\frac{C(N,K,\dist_{\X,x_0}(\cdot))-1}{\dist_{\X,x_0}(\cdot)}\m_\X. 
    \end{equation*}
\end{teo}

\subsection{The target: ${\rm CAT}(\kappa)$ spaces}

For what concerns the target space, for our harmonic map we will consider a complete ${\rm CAT}(\kappa)$ space, namely a metric space with sectional curvature bounded above by $\kappa$. Let ${\rm M}_{\kappa}$ be the \emph{model space}, namely the 2-dimensional connected, simply-connected and complete Riemannian manifold with constant sectional curvature equal to $\kappa$. Let us further denote by $\dist_\kappa$ the geodesic distance on such a space and with ${\rm D}_\kappa={\rm diam}({\rm M}_\kappa)$ its diameter, i.e.
\begin{equation*}
{\rm D}_\kappa=
    \begin{cases}
       \frac{\pi}{\sqrt{\kappa}}\qquad&{\rm if}\;\kappa>0 \\
       +\infty&{\rm if}\;\kappa\leq 0.
    \end{cases}
\end{equation*}
We also set ${\rm R}_\kappa:={\rm D}_\kappa/2$.
We have the following:
\begin{defn}[${\rm CAT}(\kappa)$ space]
   Let $(\Y,\dist_\Y)$ be a complete metric space. We say that $(\Y,\dist_\Y)$ is a ${\rm CAT}(\kappa)$ space if it is geodesic and for any triple of points $a,b,c\in\Y$ such that $\dist_\Y(a,b)+\dist_\Y(b,c)+\dist_\Y(a,c)<2{\rm D}_\kappa$ and any intermediate point $d$ between $b$ and $c$ there exist comparison points $\Bar{a},\Bar{b},\Bar{c},\Bar{d}\in {\rm M}_\kappa$ such that $\dist_\Y(a,b)=\dist_\kappa(\Bar{a},\Bar{b})$, $\dist_\Y(b,c)=\dist_\kappa(\Bar{b},\Bar{c})$, $\dist_\Y(a,c)=\dist_\kappa(\Bar{a},\Bar{c})$ and
   \begin{equation*}
       \dist_\Y(a,d)\leq\dist_\kappa(\Bar{a},\Bar{d}).
   \end{equation*}
\end{defn}

We now have a key technical Lemma holding in general ${\rm CAT}(\kappa)$ spaces which is \cite[Lemma 2.3]{ZZZ19}: we shall discuss only the case $\kappa=1$ for the sake of exposition.

\begin{lemma}
    Let $(\Y,\dist)$ be a ${\rm CAT}(1)$ space. Take any ordered sequence of points $\{P,Q,R,S\}\subset\Y$ with $\dist_\Y(P,Q)+\dist_\Y(Q,R)+\dist_\Y(R,S)+\dist_\Y(S,P)\leq 2\pi$ and let $Q_m$ be the mid-point of the geodesic joining $Q$ and $R$ (which in this case is unique). Then for any $\alpha\in[0,1]$ and $\beta>0$ we get
    \begin{align}
        \label{technical result}
        \frac{1-\alpha}{2}\biggl(4\sin^2(\dist_{QR}/2)-4\sin^2(\dist_{PS}/2)\biggr)+2\alpha\sin(\dist_{QR}/2)\biggl(2\sin(\dist_{QR}/2)-2\sin(\dist_{PS}/2)\biggr) \nonumber \\ 
        \leq \biggl[1-\frac{1-\alpha}{2}\biggl(1-\frac{1}{\beta}\biggr)\biggr]4\sin^2(\dist_{PQ}/2)+2\cos(\dist_{QR}/2)\biggl(\cos(\dist_{PQ_m})-\cos(\dist_{QQ_m})\biggr) \\  
        + \biggl[1-\frac{1-\alpha}{2}\biggl(1-\beta\biggr)\biggr]4\sin^2(\dist_{RS}/2)+2\cos(\dist_{QR}/2)\biggl(\cos(\dist_{SQ_m})-\cos(\dist_{RQ_m})\biggr). \nonumber
    \end{align}
\end{lemma}

\subsection{Sobolev spaces with metric targets and Harmonic maps}

Following \cite{GT20} (after the seminal work \cite{KS93}) we shall now introduce the Korevaar-Schoen energy and its main properties, being the main tool we need to speak about harmonic functions.

Let $u\in L^2(\Omega,\Y)$ with $\Omega\subseteq\X$ open set. We call the $2$-energy density of $u$ at scale $r$ inside $\Omega$ the quantity $ \mathbf{ks}_{2,r}[u,\Omega]:\X\to\R_+$, defined as 
\begin{equation}
    \label{p-energy density}
    \mathbf{ks}_{2,r}[u](x):=
    \begin{cases}
    \bigg(\fint_{B_r(x)}\frac{\dist_\Y^2(u(x),u(y))}{r^2}\de\m(x)\bigg)^{\frac{1}{2}}\quad&{\rm if}\;B_r(x)\subset U\\
    0&{\rm otherwise}.
    \end{cases}
\end{equation}
Moreover we introduce the \emph{total energy} of $u$ in $\Omega$ as \begin{equation}
\label{KS energy}
    {\rm E}_2[u,\Omega]:=\liminf_{r\to 0}\int_{\Omega}\mathbf{ks}_{2,r}[u,\Omega]^2(x)\de\m(x).
\end{equation}
We can now define Sobolev spaces as follows
\begin{defn}[Korevaar-Schoen space and harmonic maps]
    We say that a function $u\in L^2(\Omega,\Y)$ is in ${\rm KS}^{1,2}(\Omega,\Y)$ if ${\rm E}_2[u]<+\infty$. We say that $u$ is \emph{harmonic} in $\Omega$ if $u=\argmin_{v\in {\rm KS}^{1,2}(\Omega,\Y)} {\rm E}_{2}[v,\Omega]$.
\end{defn}
Existence of minimizers for ${\rm E}_2[\cdot,\Omega]$ has been established in the recent \cite{Sak23} (see Theorem 1.2 therein) under the condition that the image of such maps is contained in a sufficiently small ball of the target space.
%For a Lipsichitz map $u:\X\to\Y$ with source space which is a \emph{strongly rectifiable} metric measure space (which include the class of ${\rm RCD}(K,N)$ spaces) and target which is a complete metric space it is possible to speak about differential $\mathbf{md}_x(u)$ in the sense introduced by Kirchheim in \cite{Kir94} in the Euclidean setting and adapted by the authors in \cite{GT20} to the metric setting. 

We shall assume the reader to be familiar with these concepts as we are going to recall only part of \cite[Theorem 3.13]{GT20}, stating it for ${\rm RCD}$ spaces instead of the more general class of strongly rectifiable metric measure spaces.
\begin{teo}
    Let $(\X,\dist_\X)$ be an ${\rm RCD}(K,N)$ space and $(\Y,\dist_\Y)$ a complete metric space. Then for every $u\in{\rm KS}^{1,2}(\X,\Y)$ there exists a function $e_p[u]\in L^2(\X)$, called \emph{p-energy density} of $u$, such that 
    \begin{equation*}
        \mathbf{ks}_{2,r}[u]\to e_2[u]\quad\m-{\rm a.e.\;\;and\;in}\;L^2\;{\rm as}\;r\to 0.
    \end{equation*}
    In particular the $\liminf$ in \eqref{KS energy} is actually a limit.
\end{teo}
We shall now present a representation formula of the energy density $e_2[u]$ in terms of the Hilbert-Schmidt norm of the differential $|\de u|_{\rm HS}$: we will not discuss the meaning of the object $\de u$, referring to \cite{GPS18} for the details. What follows is \cite[Proposition 6.7]{GT20}.
\begin{teo}
    Let $(\X,\dist_\X,\m_\X)$ be an ${\rm RCD}(K,N)$ space and $\Omega\subset\X$ an open set. Let $(\Y,\dist_\Y)$ be a ${\rm CAT}(\kappa)$ space and $u\in{\rm KS}^{1,2}(\Omega,\Y)$, then for its energy density we have the following representation formula
    \begin{equation}
        \label{HS norm and energy density}
        e_2[u]=(d+2)^{-\frac{1}{2}}|\de u|_{\rm HS}.
    \end{equation}
\end{teo}
\begin{proof}
    Note that in \cite{GT20} the theorem is stated for $\X$ which is a strongly rectifiable space and $\Y$ which is a ${\rm CAT}(0)$ space. On one hand the proof for the case of ${\rm CAT}(\kappa)$ target is the same of the one for ${\rm CAT}(0)$ spaces, exploting the universal infinitesimal Hilbertianity of such spaces (see \cite{DMGSP18}), on the other hand we shall avoid speaking about strongly rectifiable metric measure spaces since our main results are only stated for ${\rm RCD}(K,N)$ spaces.
\end{proof}

Finally we have the following definition:
\begin{defn}[$\lambda$-convexity]
Let $(\Y,\dist_\Y)$ be a complete and geodesic metric space. We say that a function $E:\Y\to\R$ is $\lambda$-convex if for all $x,y\in\Y$ and for all geodesics $\gamma$ connecting $x=\gamma_0$ and $y=\gamma_1$ we have
\begin{equation*}
    E(\gamma_t)\leq tE(\gamma_1)+(1-t)E(\gamma_0)-\frac{\lambda}{2}t(1-t)\dist_\Y^2(\gamma_0,\gamma_1).
\end{equation*}
    
\end{defn}

\section{Main results}
\subsection{H\"older regularity of harmonic maps}
In this section we will prove H\"older regularity of our harmonic map with values in a sufficiently small ball of a ${\rm CAT}(\kappa)$ space. Note that without this assumption there may be a "big" set of discontinuity (singular set), for examples and a detailed discussion one can consult \cite{Riv95}. Since we can always renormalize the target space in such a way that it becomes a ${\rm CAT}(1)$ space, to ease the notation and the computations we shall assume $(\Y,\dist_\Y)$ to be a ${\rm CAT}(1)$ space here and in the rest of the work. 

In the following we shall prove the convexity of three functions, namely $1-\cos(\dist_{\Y,o})$, $\dist_{\Y,o}$ and $\dist^2_{\Y,o}$. The proof of the $\lambda$ convexity of the squared distance is contained \cite[Lemma3.1]{Ohta07} and the convexity of the distance $\dist_{\Y,o}$ is well-known but we shall prove them here anyway because they are natural consequences of the convexity of $1-\cos(\dist_{\Y,o})$.
\begin{prop}
\label{Convexity Cat(k)}
    Let $(\Y,\dist_\Y)$ be a ${\rm CAT}(1)$ space and consider $B_\rho(o)\subset\Y$ with $\rho<\pi/2$. Then  the distance function $\dist_{\Y,o}=\dist_\Y(o,\cdot)$ is convex on $B_\rho(o)$, $\dist^2_{\Y,o}$ is $\lambda$-convex and the function $\cos(\dist_\Y(o,\cdot))$ is $\lambda^\prime$-concave, with 
    \begin{equation*}
        \lambda = 2\cos\rho,\qquad\lambda^\prime=\cos\rho.
    \end{equation*}
    Finally $\dist_\Y(\cdot,\cdot)$ restricted to $B_{\rho/2}(o)$ is jointly convex.
\end{prop}
\begin{proof}
We show that the distance from the north pole on $\mathbb{S}^2$ is convex on the upper hemisphere. Consider three points $N, p, q\in\mathbb{S}^2$. Denote with $\dist_N(y):=\dist_{\mathbb{S}^2} (N,y)$ the distance from the north pole for every $y\in\mathbb{S}^2$ and let $\gamma$ be the geodesic connecting $p$ and $q$. By the cosine law for the sphere we can consider the triangle whose vertex are $p, q$ and $N$ and write
\begin{equation*}
    \cos(f(t))=\cos(t\dist_{\mathbb{S}^2}(p,q))\cos(\dist_N(p))+\sin(t\dist(p,q))\sin(\dist_N(p))\sin(\theta),
\end{equation*}
    where $f(t)=\dist_N(\gamma(t))$ and $\theta$ is the angle between $\gamma^\prime(0)$ and $\eta^\prime(1)$ ($\eta$ being the geodesic connecting the north pole and the point $p$). Note that we also used the fact that $\dist_{\mathbb{S}^2}(p,\gamma(t))=t\dist_{\mathbb{S}^2}(p,q)$. Now differentiate twice the previous identity to get 
    \begin{equation*}
        (\cos(f(t)))^{\prime\prime}=-\dist^2_{\mathbb{S}^2}(p,q)\cos(f(t))\leq -\dist_{\mathbb{S}^2}^2(\gamma_1,\gamma_0)\cos\rho,
    \end{equation*}
whence $\cos(f(t))$ is a $\lambda^\prime$-concave function with $\lambda^\prime=\cos(\rho)$.  Now write $f=\arccos{\cos(f)}$ and let us call $g(t):=\cos(f(t))$: we have 
\begin{equation}
\label{sombrero}
    \frac{\de^2}{\de t^2}f=\frac{(g^\prime)^2g-g^{\prime\prime}(1-g^2)}{(1-g^2)^{\frac{3}{2}}}\geq 0,
\end{equation}
meaning that $f$ is a convex function (we have used that ${\rm Im}(g)\subseteq(0,1]$ and $g^{\prime\prime}\leq 0$)- this is fully justified if $g\neq 1$, i.e. $f\neq 0$, otherwise the argument is justified by slightly moving the north pole $N$ combined with the stability properties of convexity.

For what concerns the squared distance $f^2$ just use the product rule for the derivative to get 
\begin{equation*}
    \frac{\de^2}{\de t^2}f^2= 2 |f^\prime|^2+2ff^{\prime\prime}\geq 2ff^{\prime\prime}.
\end{equation*}
Now plug \eqref{sombrero} into the previous expression to get 
\begin{equation*}
    \frac{\de^2}{\de t^2}f^2\geq 2f\bigg[\frac{(g^\prime)^2g-g^{\prime\prime}(1-g^2)}{(1-g^2)^{\frac{3}{2}}}\bigg]\geq -2f\frac{g^{\prime\prime}(1-g^2)}{(1-g^2)^{\frac{3}{2}}}\geq 2\dist^2_{\Y}(p,q)\cos\rho\frac{f}{\sin f}\geq 2\dist^2_{\mathbb{S}^2}(p,q)\cos\rho,
\end{equation*}
which is the $\lambda$ convexity with $\lambda=2\cos\rho$.

Now consider three points $x, y\in B_\rho(o)\subseteq\Y$ and let $p,q,N$ be three comparison points of $x,y,o$ in $\mathbb{S}^2$: by the ${\rm CAT}(1)$ condition we have $\dist_\Y(\Tilde{\gamma}(t),o)\leq\dist_{\mathbb{S}^2}(\gamma(t),N)$ (with ${\gamma}$ geodesic joining $p$ and $q$ and with $\Tilde{\gamma}$ geodesic joining $x$ and $y$ and), meaning that
\begin{equation*}
    \cos(\dist_\Y(\Tilde{\gamma}(t),o))\geq \cos(\dist_{\mathbb{S}^2}(\gamma(t),N)).
\end{equation*}
 The definition of comparison points together with the previous observation allows to write
\begin{equation*}
    \cos(\dist_\Y(\Tilde{\gamma}(t),o))\geq t \cos(\dist_\Y(q,o))+(1-t) \cos(\dist_\Y(p,o))+\frac{t(1-t)}{2}\dist_\Y^2(p,q)\cos\rho,
\end{equation*}
which is the sought $\lambda^\prime$-concavity with $\lambda^\prime=\cos\rho$. Analogous arguments apply for $\dist_\Y(o,\cdot)$ and $\dist_\Y^2(o,\cdot)$.

For the final part of the proof fix $x\in B_{\rho/2}(o)$ and notice that for all $y\in B_{\rho/2}(o)$ we must have $\dist_\Y(x,y)< \rho$ by triangle inequality. Therefore we can use the fact that $B_{\rho}(x)$ is convex and conclude.
\end{proof}
\begin{comment}

\begin{remark}
    Being the distance function nonnegative it is apparent that also $\frac{\dist^2_N}{2}$ is convex. 
\end{remark}
 %\bibliography{Biblio} 

Letting $R_k:=\pi/\sqrt{k}$ we have the following (see \cite{Ohta07})
\begin{teo}
\label{distance convexity}
Let $(\Y,\dist_\Y,p)$ be a pointed complete  ${\rm CAT}(\kappa)$ space and let $\eps>0$, then the function $\dist^2_o:\Bar{B}_{(R_k-\eps)/2}(o)\to[0,\infty)$ is $\lambda$-convex, namely for all $x,y\in \Bar{B}_{(R_k-\eps)/2}(o)$ and $\gamma$ geodesic between the two points we have
\begin{equation}
    \label{Lambda convexity}
    \dist^2_{o}(\gamma_t)\leq t \dist^2_{o}(y)+(1-t) \dist^2_{o}(x)-\lambda t(1-t) \dist^2(x,y)\quad\forall t\in[0,1],
\end{equation}
    where $\lambda = (R_k-\eps)\tan(\eps/2)$.
\end{teo}
\end{comment}

We recall now some lemmas of gradient flow theory on locally ${\rm CAT}(\kappa)$-spaces which will be useful to prove some Laplacian bounds. Let us start with the following, which is part of \cite[Theorem 3.3]{gn20}, to which we also refer for the relevant definitions:
\begin{teo}
    Let $\Y$ be a locally ${\rm CAT}(\kappa)$-space, $E: \Y\to \R\cup\{+\infty\}$ a $\lambda$-convex and lower semicontinuous functional. Then, the following hold:
\begin{itemize}
    \item Existence

For every $y\in D(E)$ there exists a gradient flow trajectory for $E$ starting from $y$.
\item Uniqueness and $\lambda$-contraction

For any two gradient flow trajectories $(y_t),(z_t)$ we
have
\begin{equation}
    \label{Contraction estimate}
    \dist_{\Y}(y_t,z_t)\leq e^{-\lambda(t-s)}\dist_\Y(y_s,z_s)\quad \forall t\geq s\geq 0. 
\end{equation}

\end{itemize}
\end{teo}

Then we have the following a priori estimates for the gradient flow trajectory which is \cite[Lemma 3.4]{gn20}, following the ideas contained in \cite{Pet07}:
\begin{lemma}
Let $\Y$ be locally ${\rm CAT}(\kappa)$ and $E: \Y\to \R\cup\{+\infty\}$ be a
$\lambda$-convex and lower semicontinuous functional, $\lambda\in\R$. Let $y,z\in\Y$ and consider the
gradient flow trajectories $(y_t),(z_t)$ associated with $E$.
Then, for any $t\geq s>0$, it holds
\begin{align}
    \label{A priori estimate}
    \dist^2_{\Y}(y_t,z_s)\leq&e^{-2\lambda s}\bigg(\dist^2_\Y(y,z)+2(t-s)(E(z)-E(y)) \nonumber \\
    &+2|\partial^{-}E|^2(y)\int_0^{t-s}\theta_\lambda(r)\de r-\lambda\int_0^{t-s}\dist^2_\Y(y_r,z)\de r\bigg).
\end{align}
where $\theta_\lambda(t):=\int_0^te^{-2\lambda r}\de r$.
\end{lemma}
With the previous two lemmas at hand we can prove the analogue of \cite[Lemma 4.17]{gn20} for ${\rm CAT}(\kappa)$ spaces. Below we shall denote with $\Lip_{\rm bs}(\X)$ the space of Lipschitz functions with bounded support and with $\Lip_{\rm bs}^+(\X)$ the subset of $\Lip_{\rm bs}(\X)$ made of nonnegative functions.
\begin{lemma}
\label{Gigli-Nobili CAT(k)}
    Let $(\X,\dist,\m)$ be an ${\rm RCD}(K,N)$ space, $\Y$ a locally ${\rm CAT}(\kappa)$-space and  $\Omega\subset \X$ open and bounded. Also, let $f\in\Lip(\Y)$ be $\lambda$-convex, $\lambda\in\R$, and $u\in {\rm KS}^{1,2}(\Omega,\Y)$. For $g\in\Lip_{\rm bs}(\X)^+$ define the (equivalence class of the) variation map $u_t(x):={\rm GF}^f_{tg(x)}(u(x))$ $\forall t > 0, x\in\Omega$. Then, $u_t\in {\rm KS}^{1,2}(\Omega,\Y)$ for every $t>0$ and there is a constant $C>0$ depending
on $f,g$ such that
\begin{equation}
    \label{small time behaviour}
    |\de u_t|^2_{\rm HS}\leq e^{-2\lambda t g}\bigg(|\de u|^2_{\rm HS}-2t\langle\de g, \de(f\circ u)\rangle + Ct^2\bigg)\;\m-{\rm q.o.}\,{\rm in}\,\Omega
\end{equation}
holds for every $t\in[0,1]$. In particular
\begin{equation}
    \label{limsup inequality KS}
    \limsup_{t\to 0}\frac{E^{\rm KS}(u_t)-E^{\rm KS}(u)}{t}\leq -\frac{1}{d+2}\int_{\Omega}\bigg(\lambda g|\de u|^2_{\rm HS}+\langle\de g,\de (f\circ u)\rangle\bigg)\de\m.
\end{equation}
\end{lemma}
\begin{proof}
The fact that $u_t\in L^2(\Omega,\Y)$ easily follows from the following inequalities and the fact that the support of $g$ is bounded:
\begin{align*}
    \dist_{\Y}^2(u_t(x),o)\leq2\dist_{\Y}^2(u_t(x),u(x))+2\dist_\Y^2(u(x),o) \\
    \leq 2\dist_\Y^2(u(x),o)+2te^{2|\lambda|t}\Lip^2(f)g(x),
\end{align*}
where for the second inequality we applied the a priori estimates \eqref{A priori estimate} and exploited the fact that $|\partial^{-}f|(y)\leq\Lip(f)$ for all $y\in\Y$.
Now thanks to \eqref{Contraction estimate} we have (w.l.o.g. assume $g(y)\geq g(x)$)
\begin{equation*}
    \dist^2_\Y(u_t(x),u_t(y))\leq e^{2|\lambda||g(x)-g(y)|}\dist^2_\Y\big(u(x),{\rm GF}^f_{t|g(y)-g(x)|}(u(y))\big).
\end{equation*}
Now we can use the \emph{sharp dissipation rate} of the gradient flow (see \cite[point (ii) of Theorem 3.2]{gn20}) to establish the Lipschitzianity of the map $t\to {\rm GF}^f_t(u(x))$ and get
\begin{align*}
    \dist^2_\Y\big(u(x),{\rm GF}^f_{t|g(y)-g(x)|}(u(y))\big)&\leq 2\dist^2_\Y\big(u(x),{\rm GF}^f_{t|g(y)-g(x)|}(u(x))\big) \\ 
    &+2\dist^2_\Y\big({\rm GF}^f_{t|g(y)-g(x)|}(u(x)),{\rm GF}^f_{t|g(y)-g(x)|}(u(y))\big) \\
    &\leq C_1 t^2|g(x)-g(y)|^2+2e^{2|\lambda||g(y)-g(x)|}\dist^2_\Y(u(y),u(x)) \\
    &\leq C_1 t^2 \dist^2(x,y)+ C_2\dist^2_\Y(u(y),u(x)).
\end{align*}
Dividing by $r^2:=\dist^2(x,y)$ and $\m(B_r(x))$ and integrating over $B_r(x)\subseteq\Omega$ we get
\begin{equation*}
    {\rm ks}^2_{2,r}[u_t,\Omega](x)\leq C_1 t^2 + C_2 {\rm ks}^2_{2,r}[u,\Omega](x).
\end{equation*}
The fact that $\m(\Omega)<+\infty$ allows to conclude $u_t\in {\rm KS}^2(\Omega,\Y)$. 

For what concerns estimate \eqref{small time behaviour} the proof is verbatim the one in \cite[Lemma 4.17]{gn20}.

Finally for the last point we just need to subtract from both sides of \eqref{small time behaviour} the quantity $|\de u|^2_{\rm HS}$ and then integrate over $\Omega$ and divide by $2t(d+2)$. Taking the $\limsup$ as $t\to 0^+$ and exploiting a dominated convergence argument allows to conclude with \eqref{limsup inequality KS}.
\end{proof}
The following is a generalization to ${\rm CAT}(\kappa)$ spaces of well-known inequalities holding for functions in ${\rm CAT}(0)$ spaces. We begin with the following:
\begin{prop}
\label{convex + harmonic}
   Let $(\X,\dist,\m)$ be an ${\rm RCD}(K,N)$ space and $(\Y,\dist_\Y)$ be a locally ${\rm CAT}(\kappa)$ space. Let $\Omega\subset\X$ be open and bounded and let $u:\Omega\to\Y$ be an harmonic map and $f:\Y\to\R$ be a Lipschitz and $\lambda$-convex map, then $f\circ u\in W^{1,2}(\Omega)$ and 
   \begin{equation}
       \label{Laplacian lower bound}
       \Delta(f\circ u)\geq\lambda|{\rm d}u|^2_{\rm HS}\m
   \end{equation}
   in the weak sense. In particular $\Delta(f\circ u)$ is a signed Radon measure.
\end{prop}
\begin{proof}
   The fact that $f\circ u\in W^{1,2}(\Omega)$ is well-known (see \cite{GT20}). To prove \eqref{Laplacian lower bound} first observe that being $u$ harmonic implies 
   \begin{equation*}
       \limsup_{t\to 0}\frac{E^{KS}(u_t)-E^{KS}(u)}{t}\geq 0,
   \end{equation*}
   so that \eqref{limsup inequality KS} gives
   \begin{equation*}
       \lambda\int_{\Omega}g|\de u|^2_{\rm HS}\de\m\leq - \int_{\Omega}\langle \de g,\de(f\circ u)\rangle\de\m=\int_{\Omega}\Delta(f\circ u)g\de\m
   \end{equation*}
   for all $g\in\Lip_{bs}^+(\X)$, whence \eqref{Laplacian lower bound} follows.
\end{proof}
\begin{lemma}\label{lem-local-laplacian}
    Let $(\X,\dist_\X,\m)$ be an ${\rm RCD}(K,N)$ space and $(\Y,\dist_\Y)$ a ${\rm CAT}(1)$ space. Let $u:\Omega\subset\X\to\Y$ be an harmonic mapping such that $u(\Omega)\subset B_\rho(o)$ for some $\rho<\pi/2$, then consider the function $f_o:\X\to [0,1]$ given by $f_o(x):=\cos(\dist_{\Y}(u(x),o))$. We have $f_o\in W^{1,2}(\Omega)$ and 
    \begin{equation}
    \label{local Laplacian comparison}
        \Delta f_o \leq -\cos\rho |\de u|^2_{\rm HS}
    \end{equation}
    in the weak sense in $\Omega$.
\end{lemma}
\begin{proof}
    This is indeed a consequence of Proposition \ref{Convexity Cat(k)} in combination with Proposition \ref{convex + harmonic}. Indeed one just needs to apply those results with the space $(\overline{B_\rho(o)},\dist_\Y)$, which is a ${\rm CAT}(1)$ space.
\end{proof}
We now have the following result which holds in a more general setting than the present one (see \cite[Theorem 5.4]{BiMosco95}) but we shall present it in the setting of ${\rm RCD}$ spaces to avoid further technicalities.
\begin{teo}[Elliptic Harnack inequality]
    Let $(\X,\dist,\m)$ be an ${\rm RCD}(K,N)$ space and $u:\X\to\R$ be a weakly subharmonic function in $B_{4r}(x_0)$, i.e. $u\in W^{1,2}(B_{4r}(x_0))$ and
    \begin{equation*}
        \Delta u\geq 0
    \end{equation*}
    in the weak sense in $B_{4r}(x_0)$. Then the following estimate holds
    \begin{equation}
    \label{Harnack type inequality}
        \sup_{z\in B_{r/2}(x_0)}\max\{u,0\}(z)\leq C(K^{-}r^2,N)\biggl(\frac{1}{\m(B_r(x_0))}\int_{B_r(x_0)}u^2\de\m\biggr)^{1/2},
    \end{equation}
    where $C$ is equibounded as $r\to 0^+$.
\end{teo}
\begin{remark}
    As a consequence of \eqref{Harnack type inequality} we get that any weakly subharmonic function is locally bounded from above.
\end{remark}

We shall now introduce the following notation: for a function $v:\X\to\R$ we set 
\begin{equation*}
    v_R:=\fint_{B_R(x_0)}v\de\m,
\end{equation*}
where $x_0\in\X$ is a point which will be clear from the context. We further set 
\begin{equation*}
    v_{+,R}:=\sup_{x\in B_R(x_0)}\max\{v,0\}(x)
\end{equation*}
The following is a combination of \cite[Corollary 1]{Jost97} and \cite[Lemma 7]{Jost97}:
\begin{cor}
\label{Corollary Jost}
    Let $u:\X\to\R$ be as in the previous Theorem and nonnegative, then there exists $\delta_0>0$ independent of $R$ such that 
    \begin{equation*}
        \sup_{B_{R}(x_0)} u\leq (1-\delta_0)u_{+,4R} + \delta_0u_R.
    \end{equation*}
   Moreover if $\eps\in(0,1/4)$ there exists $m\in\N$ (independent of $u$ and $\eps$) such that 
   \begin{equation}
       \label{Corollary Jost II}
       u_{+,\eps^mR}\leq\eps^2u_{+,R}+(1-\eps^2)u_{R^\prime}
   \end{equation}
   where $R^\prime$ (possibly depending on $\eps$ and $u$) is such that  $\eps^mR\leq R^\prime\leq R/4$.
\end{cor}

We proceed recalling another useful lemma which again extends to the context of ${\rm CAT}(\kappa)$ spaces without modifications:
\begin{lemma}
    Let $(\X,\dist,\m)$ be an ${\rm RCD}(K,N)$ space and let $(\Y,\dist_\Y,o)$ be a pointed complete metric space, then for every $u\in {\rm KS}^{1,2}(\X,\Y)$ there exists $C=C({\rm diam}(\Omega),K,N)\geq 1$ such that for every $r>0$ and $p\in\Omega$ for which $B_{rC}(p)\subseteq\Omega$ we have 
    \begin{equation}
        \label{Poincarè inequality}
        \fint_{B_r(p)}\int_{B_r(p)}\dist^2_\Y(u(x),u(y))\de\m(x)\de\m(y)\leq Cr^2\int_{ B_{rC}(p)}e_2^2[u]\de\m.
    \end{equation}
\end{lemma}
\begin{proof}
  The proof can be found in \cite[Lemma 4.9]{Gu017v5}. 
\end{proof}
The next Lemma is basically \cite[Lemma 8]{Jost97} adapted to ${\rm CAT}(\kappa)$ setting.

\begin{lemma}
\label{Technical Lemma}
    Let $(\X,\dist,\m)$ be an ${\rm RCD}(K,N)$ space, $\Omega\in\X$ an open set, and $(\Y,\dist_\Y)$ be a ${\rm CAT}(1)$ space. Let $u:\Omega\to\Y$ be an harmonic map with values in $B_\rho(o)$ with $\rho<\pi/2$ and let $B_{4R}(x_0)\subset\subset\Omega$, then 
    \begin{equation*}
        R^2\fint_{B_{R}(x_0)}|\de u|^2_{\rm HS}\de\m\leq C(v_{+,4R}-v_{+,R}),
    \end{equation*}
    where $v(x)=\dist^2_\Y(u(x),o)$ and $C=C({\rm diam}(\Omega),K,N)$.
\end{lemma}
\begin{proof}
    To begin with let us consider a mollified version of the Green function (whose existence can be proved for instance via Lax-Milgram theorem) which solves in the weak sense the following
     \begin{equation*}
    \begin{cases}
       
            -\Delta G_{p}=\frac{\chi_{B_R(p)}}{\m(B_R(p))}\qquad{\rm on}\;B_{2R}(p) \\
            
            G_p=0\qquad{\rm on}\;B^c_{2R}(p).
    \end{cases}
      \end{equation*} 
      We have (we shall omit the point $p$ center of the ball)
      \begin{equation}
          \label{Approximate Green function}
          \int_{B_{2R}}\langle\de\varphi,\de G_p\rangle\de\m = \fint_{B_R}\varphi\de\m
      \end{equation}
      for all $\varphi\in\Lip_{bs}(\X)$ with ${\rm supp}\varphi\subset\subset B_{2R}(p)$. Now following \cite[Section 6]{BiMosco95} we define for convenience a rescaled version of $G$, namely we set
      \begin{equation*}
          G_{p,R}:=\frac{\m(B_R(p))}{R^2}G_p,
      \end{equation*}
      which satisfies
      \begin{equation*}
           \int_{B_{2R}}\langle\de\varphi,\de G_{p,R}\rangle\de\m = \frac{1}{R^2}\int_{B_R}\varphi\de\m
      \end{equation*}
      and the following estimates (again we refer to \cite[Theorem 6.1]{BiMosco95}, which deals with more general metric spaces which  include the class of ${\rm RCD}(K,N)$ spaces)
      \begin{align*}
          0<C_1\leq G_{p,R}\quad{\rm on}\;B_R, \\
          0\leq G_{p,R}\leq C_2\quad{\rm on}\;B_{2R},
      \end{align*}
      where $C_1,C_2$ only depend on $K$,$N$ and ${\rm diam}(\Omega)$.
      Now we can define $z:=v-v_{+,4R}$ and write, exploiting \eqref{Laplacian lower bound} for $f(\cdot)$ equal to  $\dist_\Y^2(\cdot,o)$ with $\lambda=2\cos\rho$ by Proposition \ref{Convexity Cat(k)},
      \begin{equation*}
          \lambda\int_{B_{2R}}|\de u|^2_{\rm HS}G_{p,R}^2\de\m\leq\int_{B_{2R}}(\Delta z) G^2_{p,R}\de\m=-2\int_{B_{2R}}\langle\de z,\de G_{p,R}\rangle G_{p,R}\de\m.
      \end{equation*}
      Now we can use the Leibniz rule for the differential $\de(G_{p,R}v)=G_{p,R}\de z+z\de G_{p,R}$ and write
      \begin{equation*}
          \lambda\int_{B_{2R}}|\de u|^2_{\rm HS}G_{p,R}^2\de\m\leq-2\int_{B_{2R}}\langle\de G_{p,R},\de(G_{p,R}z)\rangle\de\m+2\int_{B_{2R}}\langle\de G_{p,R},\de G_{p,R}\rangle z\de\m.
      \end{equation*}
      Being $z\leq 0$ we can neglect the second term and obtain
      \begin{align*}
            \lambda\int_{B_{2R}}|\de u|^2_{\rm HS}G_{p,R}^2\de\m&\leq-2\int_{B_{2R}}\langle\de G_{p,R},\de(G_{p,R}z)\rangle\de\m=-\frac{1}{R^2}\int_{B_R}G_{p,R}z\de\m \\
            &\leq -\frac{C_1\m(B_R)}{R^2}(v_R-v_{+,4R})=\frac{C_1\m(B_R)}{R^2}(v_{+,4R}-v_R)
      \end{align*}
      where we used the definition of the mollified Green function. Finally, applying Corollary \ref{Corollary Jost}, we get the thesis.
\end{proof}
We are now in position to prove the desired H\"older continuity of harmonic maps.
\begin{teo}
\label{Holder continuity}
    Let $u:\Omega\subseteq\X\to\Y$ be an harmonic map such that ${\rm Im}(u)\subseteq B_{\rho}(o)$ with $\rho<\pi/2$ and with $(\X,\dist,\m)$ which is an ${\rm RCD}(K,N)$ space and $\Y$ which is a  ${\rm CAT}(1)$ space. Then $u$ is locally H\"older continuous in $\Omega$.
\end{teo}
\begin{proof}
    The proof closely follows \cite[Theorem]{Jost97}. Let us fix $x_0\in\Omega$ in such a way that $B_{4R}(x_0)\subset\subset\Omega$. Let us define the mean of $u$ on a ball centered at $x_0$ with radius $r$, denoted by $\Bar{u}_r$, as one of the minimums of
    \begin{equation*}
        \Y\ni q\mapsto\fint_{B_r(x_0)}\dist_\Y(u(x),q)\de\m(x).
    \end{equation*}
Finally set $v_p(x):=\dist_\Y^2(u(x),p)$ where $p\in\Y$ will be chosen later and $w(x):=\dist^2_{\Y}(u(x),\Bar{u}_{R/4})$. We want to exploit the result in Corollary \ref{Corollary Jost}: let us therefore fix $\eps\leq 1/10$ so that $\eps^mR\leq R^\prime\leq R/4$ and estimate as follows 
\begin{equation*}
    w_{R^\prime}^m=\frac{1}{\m(B_{R^\prime}(x_0))}\int_{B_{R^\prime}(x_0))}\dist^2_{\Y}(u(x),\Bar{u}_{R/4})\de\m(x)\leq \frac{C}{\m(B_{R/4}(x_0)}\int_{B_{R/4}(x_0)}\dist^2_{\Y}(u(x),\Bar{u}_{R/4})\de\m(x)
\end{equation*}
    where $C$ is independent of $R$, exploiting the (uniformly) doubling property of the measure $\m$ on $\Omega$.
    Now applying Poincaré inequality to the previous expression we get 
\begin{equation*}
     \frac{C}{\m(B_{R/4}(x_0))}\int_{B_{R/4}(x_0)}\dist^2_{\Y}(u(x),\Bar{u}_{R/4})\de\m(x)\leq C_1\frac{R^2}{\m(B_R(x_0))}\int_{B_{R/(4\lambda)}(x_0)}|\de u|_{\rm HS}^2\de\m,
\end{equation*}
for some $\lambda\in (0,1)$. Now we shall apply Lemma \ref{Technical Lemma} and the doubling inequality again to obtain 
\begin{equation}
\label{key estimate holder continuity}
   w^m_{R^\prime}\leq C\big(v_{p,+,R/\lambda}-v_{p,+,R/4\lambda}\big). 
\end{equation}

Choose now $p\in {\rm conv}\big(u(B_{\eps^mR}(x_0)\big)$ so that we have 
\begin{equation*}
    \sup_{x\in B_{\eps^mR}(x_0)}\dist^2_{\Y}(u(x),p)\leq 2 \sup_{x\in B_{\eps^mR}(x_0)}\dist^2_{\Y}(u(x),\Bar{u}_{R/4})+2\dist^2_{\Y}(\Bar{u}_{R/4},p)\leq 4 \sup_{x\in B_{\eps^mR}(x_0)}\dist^2_{\Y}(u(x),\Bar{u}_{R/4})
\end{equation*}
and at the same time 
\begin{equation*}
     \sup_{x\in B_{R}(x_0)}\dist^2_{\Y}(u(x),\Bar{u}_{R/4})\leq 4 \sup_{x\in B_{R}(x_0)}\dist^2_{\Y}(u(x),p)
\end{equation*}
Combining estimate \eqref{key estimate holder continuity} and the result of Corollary \ref{Corollary Jost} we get 
\begin{align*}
     \sup_{x\in B_{\eps^mR}(x_0)}\dist^2_{\Y}(u(x),\Bar{u}_{R/4})
     \leq 4\eps^2\sup_{B_{R}(x_0)}\dist^2_\Y(u(x),\Bar{u}_{R/4}) + C\big(v_{p,+,R/\lambda}-v_{p,+,R/4\lambda}\big) \\ 
     \leq 16\eps^2\sup_{x\in B_{R}(x_0)}\dist^2_{\Y}(u(x),p)+C\big(v_{p,+,R/\lambda}-v_{p,+,\eps^mR}\big),
\end{align*}
where in the last line we also used that $\eps^m\leq (1/8)^m\leq 1/4\leq 1/4\lambda$. In the end we obtain
\begin{equation*}
      \sup_{x\in B_{\eps^mR}(x_0)}\dist^2_{\Y}(u(x),p)\leq 64\eps^2\sup_{x\in B_{R}(x_0)}\dist^2_{\Y}(u(x),p)+C\big(v_{p,+,R/\lambda}-v_{p,+,\eps^mR}\big).
\end{equation*}

Setting $ \omega(r):=\sup_{x\in B_r(x_0)}\dist^2_\Y(u(x),p)$ we can rewrite the previous inequality as 
\begin{equation*}
    (1+C)\omega(\eps^mR)\leq 64\eps^2\omega(R)+C\omega(R/\lambda)\leq (64/100+C)\omega(R/\lambda),
\end{equation*}
which means 
\begin{equation*}
    \omega(\eps^mR)\leq c\omega(R/\lambda),
\end{equation*}
where $\eps$ and $\lambda$ are fixed and $c<1$. By an iteration of the latter estimate (holding for every $R\leq R_0$ for which $B_{R_0}(x_0)\subset\subset\Omega$) we get
\begin{equation*}
    \frac{\omega(r)}{r^\alpha}\leq C\frac{\omega(R_0)}{R_0^\alpha},
\end{equation*}
where $\alpha\in (0,1)$, $C>0$ and $r\leq R_0$. Choosing $p=\Bar{u}_r$ we get 
\begin{equation*}
    \sqrt{\omega(r)}\leq{\rm osc}(u,B_r(x_0)) \leq 2  \sqrt{\omega(r)}
\end{equation*}
and this proves the (local) H\"older continuity of $u$.
\end{proof}

\subsection{Higher integrability  of energy densities}

Let $\Omega\subset X$ be an open bounded set in an ${\rm RCD}(K,N)$ space with $X\setminus \Omega\not=\emptyset$ , $K\in\mathbb R$ and $N\in[1,\infty)$.  Let $(Y,{\dist}_{\Y})$ be a ${\rm CAT}(1\kappa)$ space with $\kappa>0$ and suppose that $u\in{\rm KS}(\Omega;Y)$ is an harmonic map with values in a ball $B_\rho(o)\subset Y$ with $\rho\in(0,\frac{\pi}{2\sqrt{\kappa}}).$ We shall always fix a H\"older continuous representative of $u$.

Let us recall some notations in \cite{AndersonTH17}.
\begin{defn}
    For any $q>1$, a nonnegative $\m$-measurable function $w$ on $\Omega$ belongs to the {\it weak $q$-Reverse H\"older class} $RH^{weak}_q$ if there exists a constant $C_q$ such that 
    \begin{equation*}
        \left(\fint_B w^q\de\m\right)^{1/q}\leqslant  C_q \int_{2B}w\de\m
    \end{equation*} 
for all ball $B:=B_r(y)$ with $2B:=B_{2r}(y)\subset \Omega$.
\end{defn}

We need the following Gehring lemma, see \cite[Propostion 6.2]{AndersonTH17} and \cite[Theorem 3.1]{Maasalo07}.
\begin{lemma}\label{gehring lemma}
    If $1<q<\infty$ and $w\in RH^{weak}_q$, then there exists $\varepsilon>0$ such that $w\in RH^{weak}_{q+\varepsilon}.$    \end{lemma}
Now we will prove the higher integrability of energy density.
\begin{teo}
    \label{higher-integrable}
    Let $\Omega, Y$ and $u$ be as above. Then there exists an $\varepsilon=\varepsilon(N,K,{\rm diam}(\Omega),\rho)>0$ such that $|\de u|_{\rm HS}\in W^{1,2+\epsilon}_{\rm loc}(\Omega)$ and  
    \begin{equation}\label{doubling-energy}
    \left(\fint_{B} |\de u|^{2+\varepsilon}_{\rm HS}\de\m\right)^{\frac{2}{2+\varepsilon}}\leqslant C_\varepsilon\fint_B|\de u|_{\rm HS}^2\de\m
    \end{equation}
    for any ball $B$ with $2B\subset \Omega,$ where the constant $C_\epsilon>0$ depends only on $\epsilon$.
    \end{teo}
    
\begin{proof}
    Fix any ball $B$ with $2B\subset\Omega$, then by Lemma \ref{lem-local-laplacian}, we have 
    $$\Delta (f_o-a)\leqslant -\cos\rho |\de u|_{\rm HS}^2,\quad \forall a\in\mathbb R,$$
where $f_o(x)=\cos({\dist}_{\Y}(u(x),o))$. Let $\phi:\Omega\to [0,1]$ be a cut-off function with $\phi=1$ on $B$, $\phi=0$ out of $\frac{3}{2}B$, and $$|\nabla \phi|\leqslant C_1r^{-1},\quad |\Delta \phi|\leqslant C_2r^{-2},$$
where the constants $C_1,C_2$ depend only on $K,N$ and ${\rm diam}(\Omega)$. Then we get
$$\int_B|\de u|^2_{\rm HS}\de\m \leqslant\int_{\frac{3}{2}B}|\de u|^2_{\rm HS}\phi \de\m \leqslant \frac{C_{3}}{r^2}\int_{\frac 3 2 B} |f_o-a|\de\m,$$
for all $a\in\mathbb R$, where $C_3=C_2/\cos\rho.$ It is well-known that a weak $(1,2)$-Poincar\'e inequality holds on ${\rm RCD}(K,N)$ spaces and since the weak $(1,s)$-Poincar\'e inequality is an open ended condition (see \cite[Theorem 1.0.1]{KeithZhong2008}), there exists  a number $s_0\in(1,2)$ such that the weak $(1,s_0)$-Poincar\'e inequality holds on  ${\rm RCD}(K,N)$ spaces. Therefore, we have 
$$\inf_{a\in\mathbb R}\fint_{\frac 3 2 B}|f_o-a|\de\m\leqslant C_{K,N,{\rm diam}(\Omega),s_0}\cdot r\left(\fint_{2B} |\nabla f_o|^{s_0}\de\m\right)^{1/s_0}.$$
Combining the above two inequalities, we conclude that 
$$\left(\fint_B|\de u|^2_{\rm HS}\de\m \right)^{1/2}\leqslant C_4\left(\fint_{2B}|\nabla f_o |^{s_0}\de\m\right)^{1/s_0}\leqslant C_4 \left(\fint_{2B}|\de u |_{\rm HS}^{s_0}\de\m\right)^{1/s_0},$$
where we have used $|\nabla f_o|\leqslant |\sin \dist_Y(o,u)|\cdot |\nabla \dist_Y(o,u)|\leqslant   |\de u|_{\rm HS}.$ Now, applying \ref{gehring lemma} to $|\de u|^{s_0}_{\rm HS}$, we obtain
$|\de u|_{\rm HS}$ is in $W^{2+\varepsilon}_{\rm loc}(\Omega)$, and moreover 
$$\left(\fint_B|\de u|^{2+\varepsilon}_{\rm HS}\de\m \right)^{1/(2+\varepsilon)}\leqslant    C_\varepsilon \left(\fint_{2B}|\de u |_{\rm HS}^{s_0}\de\m\right)^{1/s_0}\leqslant C_\varepsilon \left(\fint_{2B}|\de u |_{\rm HS}^{2}\de\m\right)^{1/2},$$
since $s_0<2$ and H\"older inequality.
\end{proof}

\subsection{Auxiliary results}

In this section we shall work under the following assumptions:
    \begin{enumerate}
        \item $(\X,\dist_\X,\m)$ is an ${\rm RCD}(K,N)$ space with essential dimension $d\in\N$. 
        \item $\Omega\subset\X$ is an open bounded set with $X\setminus\Omega\not=\emptyset$, which is equivalent to  $\m(\X\setminus\Omega)>0$. 
        \item $(\Y,\dist_\Y)$ is a ${\rm CAT}(1)$ space: the results obtained for general ${\rm CAT}(\kappa)$ spaces will be obtained by a rescaling of the distance function. 
        \item $u\in{\rm KS}^{1,2}(\Omega;\Y)$ is harmonic with values in a ball $B_\rho(o)\subset\Y$ with $\rho<\frac{\pi}{2}$. Finally we shall fix Borel representatives of $u$ (the H\"older continuous one) and of $e_2[u]$ (and of $|\de u|_{\rm HS}$).
        \item Let $\Omega^\prime\subset\subset\Omega$ be open and consider $r>0$ and $\hat{x}\in\Omega^\prime$ such that $B_{4r}(\hat{x})\subset\Omega^\prime$ and $\|u\|_{C^\alpha(\Omega^\prime)}r^\alpha<\pi/10$. Finally call $B=B_r(\hat{x}),$ $2B:=B_{2r}(\hat{x})$ and $B^\prime=B_{3r/2}(\hat{x})$.

    \end{enumerate}

Let us first define $F:\R\to\R$ as the following
\begin{equation*}
    F(t):=2\sin\bigg(\frac{t}{2}\bigg)+4\sin^2\bigg(\frac{t}{2}\bigg)
\end{equation*}
and observe that $F$ is such that $F^\prime, F^{\prime\prime}\geq 0$ on $[0,\pi/2]$. With a little abuse of notation let us also set 
\begin{equation*}
    F(z,w):=2\sin\bigg(\frac{\dist_\Y(z,w)}{2}\bigg)+4\sin^2\bigg(\frac{\dist_\Y(z,w)}{2}\bigg)
\end{equation*}
for any $z,w\in\Y$.

 We introduce the following quantities, in order to produce an Hopf-Lax formula for the function $u$,
\begin{equation*}
\label{Modified Hopf-Lax}
    f(x,y)= 
    \begin{cases}
       -F(u(x),u(y))\quad&{\rm if}\,x,y\in B^\prime \\
       -6\quad&{\rm otherwise}.
    \end{cases}
\end{equation*}
Notice that $f$ is lower semiconinuous since $F$ is bounded between $0$ and $6$. We call $f_t$ the $p$-Hopf-Lax semigroup applied to the function $f$, namely we set
\begin{equation}
\label{Hopf-Lax of distance function}
    f_t(x):=\inf_{y\in\X}\bigg[\frac{\dist_{\X}^p(x,y)}{pt^{p-1}}+f(x,y)\bigg],
\end{equation}
where we avoid to include $p$ in the definition of $f_t$ to lighten the notation.
Notice that $0\geq f_t(x)\geq -6$ for every $x\in\X$. Moreover the infimum in \eqref{Hopf-Lax of distance function} is actually a minimum (this follows by Weierstrass theorem exploiting the semicontinuity of the function we are minimizing). We also have a quantitative estimate for where to find a minimum, indeed denoting with $y_{t,x}$ a minimizer for $f_t(x)$, choosing $x$ as a competitor, we get
\begin{equation*}
    f_t(x)\leq \frac{\dist_\X^p(x,y_{t,x})}{pt^{p-1}}+f(x,y_{t,x})\leq 0.
\end{equation*}
This means $\dist_\X(x,y_{t,x})\leq (6pt)^\frac{p-1}{p}$ so that there exists $t_*=t_*(p)>0$ such that we have
\begin{equation*}
      f_t(x):=\inf_{y\in B_{12\sqrt{t}}(x)}\bigg[\frac{\dist_\X^p(x,y)}{pt^{p-1}}-F(u(x),u(y))\bigg]\qquad\forall x\in B
\end{equation*}
for $t\in (0,t_*)$.

Now set 
\begin{equation*}
    S_t(x):=\bigg\{y\in \X:\;f_t(x)=\frac{\dist_\X^p(x,y)}{pt^{p-1}}-F(u(x),u(y))\bigg\}
\end{equation*}
and observe that the latter set is non-empty if $t<t_*$. Finally set 
\begin{equation*}
    L_t(x):=\min_{y\in S_t(x)}\dist_\X(x,y)\qquad{\rm and}\qquad D_t(x):=\frac{L_t^p(x)}{pt^{p-1}}-f_t(x).
\end{equation*}
We now present a slight modification of \cite[Lemma 4.1]{ZZZ19} since we still don't know that the map $u$ is Lipschitz continuous but we have H\"older regularity instead: if the map is assumed to be Lipschitz the proof works in the same way replacing $\alpha$ with $1$.
\begin{lemma}
\label{Holder and Lipschitz time estimates}
    With the above notation and assumptions, we have
    \begin{enumerate}
        \item $f_t$ is H\"older continuous on $B$. 
        \item $L_t$ and $D_t$ are lower semicontinuous. 
        \item There exists a constant $C=C(p,\|u\|_{C^\alpha},k)>0$ such that 
        \begin{equation}
            \label{Modified Lemma}
            L_t\leq C t^{\beta},\qquad D_t\leq \Tilde{C}t^{\beta^\prime},\qquad -f_t\leq \Tilde{C}t^{\beta^\prime}\quad{\rm on}\; B,
        \end{equation}
    \end{enumerate}
where $\beta=(p-1)/(p-\alpha)$, $\beta'=\alpha \beta$, and the constant $C$  depends on $p$, the H\"older norm of $u$, $\|u\|_{C^\alpha}$.  
\end{lemma}
\begin{proof}
     The proof of $(1)$ is immediate since the infimum of equi-H\"older functions is H\"older.

    The proof of $(2)$ is contained in \cite[Lemma 4.1]{ZZZ19}.

    For the proof of $(3)$ consider $y_t(x)\in S_t(x)$ such that $L_t(x)=\dist_\X(x,y_t)$. We get, using $\sin\theta\leq\theta$ for $\theta>0$ and that $\dist_\Y(u(x),u(y_t))\leq \pi$,
    \begin{align*}
        D_t(x)=\frac{L_t^p(x)}{pt^{p-1}}-f_t(x)=F(x,y_t)\leq \dist_\Y(u(x),u(y_t))+\dist_\Y^2(u(x),u(y_t)) \\ 
        \leq (1+\pi)\dist_\Y(u(x),u(y_t))\leq(1+\pi)\|u\|_{C^\alpha}L_t^\alpha(x)=C L_t^\alpha(x).
    \end{align*}
    At the same time, being $f_t\leq 0$ we have 
    \begin{equation*}
        \frac{L_t^p(x)}{pt^{p-1}}\leq D_t(x)\leq C L_t^\alpha(x)
    \end{equation*}
    so that we get $L_t\leq C t^\beta$, with $\beta=(p-1)/(p-\alpha)$. For $D_t$ we have instead
    \begin{equation*}
        D_t(x)\leq C L_t^\alpha(x)\leq \Tilde{C}t^{\beta^\prime}
    \end{equation*}
    with $\beta^\prime=\alpha\beta$. Finally for $f_t$ we have, since $-f_t\leq D_t$,
    \begin{equation*}
        -f_t\leq\Tilde{C}t^{\beta^\prime}.
    \end{equation*}
\end{proof}

To establish the key variational inequality we shall exploit the following simple but useful lemma
\begin{lemma}
\label{Distributional laplacian bound for f}
    With the above assumptions we have 
    \begin{equation*}
        \Delta f(\cdot,y)\leq 0\quad{\rm on}\;B
    \end{equation*}
    in the weak sense, for all $y\in B$.
\end{lemma}
\begin{proof}
    Thanks to the assumptions it is sufficient to compute $\Delta F(u(\cdot),u(y))$ in the weak sense. By the chain rule we get 
    \begin{equation*}
        \Delta F(u(\cdot),u(y))=F^{\prime\prime}|\nabla\dist_\Y(u(\cdot),u(y))|^2+F^\prime\Delta\dist_\Y(u(\cdot),u(y)),
    \end{equation*}
    whence the claim follows by the nonnegativity of the factors on the right hand side (recall that the Laplacian of $x\mapsto\dist_\Y(u(x),u(y))$ is nonnegative thanks to Proposition \ref{convex + harmonic}).
\end{proof}

We now have a lemma on the heat flow Laplacian of the Hopf-Lax semigroup (the idea is from \cite{MS21}, see also \cite{Gig23} and \cite{MS22})
\begin{lemma}
    Let $f:\X\to\R$ be a bounded Borel function. Assume that for some $x,y\in\X$ we have
    \begin{equation}
    \label{equality in Hopf-Lax}
        Q_t^p f(x)=f(y)+\frac{\dist^p_\X(x,y)}{pt^{p-1}}.
    \end{equation}
    Then
    \begin{equation}
        \label{Heat flow bound Hopf-Lax}
        \Delta Q_t^pf(x)\leq \Delta f(y)-K\frac{\dist^p_\X(x,y)}{t^{p-1}}.
    \end{equation}
    holds in the heat flow sense.
\end{lemma}
\begin{proof}
    First of all let $\pi_s\in\mathcal{P}(\X\times\X)$ be an optimal transport plan between $h_s\delta_x\in\mathcal{P}(\X)$ and $h_s\delta_y\in\mathcal{P}(\X)$ for the cost $\dist_\X^p$.
    Moreover we have the following estimate, which is the Wasserstein contractivity of the heat flow (holding in general ${\rm RCD}(K,\infty)$ spaces, see \cite{AmbrosioGigliSavare12}),
    \begin{equation}
        \label{Wasserstein contractivity}
        W^p_p(h_s\delta_x,\delta_y)\leq e^{-pKs}\dist_\X^p(x,y).
    \end{equation}
    We can now estimate as follows
    \begin{align*}
        h_sQ_tf(x)&=\int_{\X}Q_tf(z)\de h_s\delta_x(z)=\int_{\X\times\X}Q_tf(z)\de\pi_s(z,z^\prime) \\
        &\leq\int_{\X\times\X}\bigg[f(z^\prime)+\frac{\dist^p(z,z^\prime)}{pt^{p-1}}\bigg]\de\pi_s(z,z^\prime) \\
        \text{(by optimality of $\pi_s$)}\qquad\qquad&=\int_\X f(z^\prime)\de h_s\delta_y(z)+\frac{1}{pt^{p-1}}W^p_p(h_s\delta_x,h_s\delta_y) \\
        &=h_sf(y)+\frac{1}{pt^{p-1}}W^p_p(h_s\delta_x,h_s\delta_y).
    \end{align*}
   Finally applying \eqref{Wasserstein contractivity} to the previous inequality we get (note that the following would hold for any $w\in\X$ in place of $y$)
   \begin{equation}
       \label{inequality case Hopf-Lax}
         h_sQ_tf(x)\leq h_sf(y)+\frac{e^{-pKs}}{pt^{p-1}}\dist_{\X}^p(x,y).
   \end{equation}
   Subtracting \eqref{equality in Hopf-Lax} from \eqref{inequality case Hopf-Lax}, dividing by $s>0$ and taking the $\limsup$ as $s\to 0$ finally gives \eqref{Heat flow bound Hopf-Lax}.
\end{proof}

We now proceed with a refinement of \eqref{local Laplacian comparison}, following \cite[Proposition 1.17]{Ser95}, which will be crucial for obtaining an elliptic inequality involving the function $f_t$.

\begin{lemma}
\label{Key Lemma}
Let $u:\Omega\to\Y$ be an harmonic map with $\Omega\subset\X$ open set, $(\Y,\dist_\Y)$ which is a ${\rm CAT}(1)$ space and ${\rm Im}(u)\subseteq B_\rho(o)$ with $o\in\Y$, $\rho<\pi/2$. Let further $f_o(x):=\cos(\dist_{\Y}(u(x),o))$, then we have $f_o\in W^{1,2}(\Omega)$ and 
\begin{equation}
    \label{key inequality}
    \Delta f_o\leq -f_o |\de u|_{\rm HS}^2=-f_o(n+2) e_2^2[u]\qquad{\rm in}\;\Omega
\end{equation}
    in the weak sense.
\end{lemma}

\begin{proof}
 Let us first set $R(x):=\dist_{\Y}(u(x),o)$, denote with $x\to G_{t}^{u(x),o}$ the map which associates to each $x\in\Omega$ the point at time $t$ lying in the geodesic (recall that geodesics are unique in our case) connecting $o$ and $u(x)$. Finally set $u_\eta:=G_{\eta}^{u,o}$ where $\eta\in W^{1,2}(\Omega)\cap C_c(\Omega)$ is such that $0\leq\eta\leq 1$: then by  \cite[Lemma 3.8]{Sak23} we have
\begin{equation}
    \label{Lemma Sakurai}
    e_2^2[u_{\eta t}]\leq \frac{\sin^2\big[(1-\eta t)R\big]}{\sin^2R}(e_2^2[u]-e_2^2[R])+e_2^2[(1-\eta t)R]
\end{equation}
 $\m$-a.e. in $\Omega$, where $t$ is a positive parameter that we will eventually send to zero. Now we shall use the duplication formula for the sinus to get 
\begin{equation*}
     |\de u_{\eta t}|_{\rm HS}^2\leq \bigg[\cos^2(t\eta R)+\frac{\sin^2(t\eta R)\cos^2R}{\sin^2R}-\frac{\cos R\sin(2t\eta R)}{\sin R}\bigg](|\de u|^2_{\rm HS}-|\de R|^2_{\rm HS})+|\de R-t\eta\de R|^2_{\rm HS}.
\end{equation*}
Note that we have simultaneously used that $|\de u|^2_{\rm HS}=(n+2)e_2^2[u]$ (recall that if $f:\X\to\R$ then $|\de f|=|\de f|_{\rm HS}$).
We proceed integrating over $\Omega$, we divide by $t$ and exploit the fact that $E_2(u_{t\eta})-E_2(u)\geq 0$ (as $u$ is harmonic) together with the asymptotics of the involved functions to get 
\begin{equation*}
    0\leq \int_\Omega\bigg[-\eta R\frac{\cos R}{\sin R}|\de u|^2_{\rm HS}+\eta R\frac{\cos R}{\sin R}|\de R|^2-\langle\de R,\de(\eta R)\rangle\bigg]\de\m.
\end{equation*}
We can now use the following identity
\begin{equation*}
    \big\langle\nabla\bigg(\eta\frac{R}{\sin R}\bigg),\nabla\cos R\big\rangle = \eta R\frac{\cos R}{\sin R}|\de R|^2-\langle\de R,\de(\eta R)\rangle
\end{equation*}
to get 
\begin{equation*}
    0\leq\int_{\Omega}-\eta R\frac{\cos R}{\sin R}|\de u|^2_{\rm HS}+\big\langle\nabla\bigg(\eta\frac{R}{\sin R}\bigg),\nabla\cos R\big\rangle\de\m.
\end{equation*}
    Note that now we can choose the magnitude of $\eta$ to be whatever we want since the inequality doesn't change if we divide everything by a positive constant. Now pick $\varphi\in\Lip_c(\Omega)$ nonnegative and set $\eta:=\varphi R/\sin R$: it is clear that $\eta\in W^{1,2}(\Omega)\cap C_c(\Omega)$ because it is the product of a bounded $W^{1,2}(\Omega)$ and continuous function and a Lipschitz function with compact support. Finally this means that for all $\varphi\in \Lip_c(\Omega)$ nonnegative we have
    \begin{equation*}
        \int_\Omega\varphi\cos R|\de u|^2_{\rm HS}\de\m\leq\int_{\Omega}\big\langle\nabla\varphi,\nabla\cos R\big\rangle\de\m.
    \end{equation*}
    The latter is the conclusion.
\end{proof}

Finally define some parametric functions depending on the distance of the target space $\dist_\Y$ and deduce some Laplacian bounds on them that we shall exploit later in the proof of the "good" distributional bound.
\begin{lemma}
\label{Key Lemma special function}
    Let $u:\Omega\subset\X\to\Y$ be an harmonic map with ${\rm Im}(u)\subset B_\rho(o)$ and $\rho<\pi/2$. Consider for any $z\in\Omega$ and $y\in\Y$  the function
    \begin{equation*}
        w_{a,b,y,z}(x)=a\dist_\Y^2(u(x),u(z))+b\cos(\dist_\Y(u(x),y)).
    \end{equation*}
    For $\m-$a.e. $x_0\in\Omega$ we have 
    \begin{align*}
       \Delta w_{a,b,o,x_0}(x_0)&\leq \big(2a-b\cos(\dist_\Y(u(x_0),o))\big)(n+2)e_2^2[u](x_0) \\
        &=\big(2a-b\cos(\dist_\Y(u(x_0),o))\big)|\de u|^2_{\rm HS}(x_0)
    \end{align*}
    in the heat flow sense.
\end{lemma}
\begin{proof}
        First of all we shall notice that \cite[Proposition 3.3]{MS22} holds also in this setting with the same proof since by Lemma \ref{Holder continuity} we have the (H\"older) continuity of $u$. Therefore we have 
        \begin{equation}
        \label{First part}
            h_t(\dist_\Y^2(u(\cdot),u(x_o))(x_0)=2|\de u|^2_{\rm HS}(x_0)t+o(t)\qquad{\rm as}\; t\to 0^+.
        \end{equation}
        for $\m$-a.e. $x_0$.
        Secondly by the results contained in \cite{Gigli_2023} and Lemma \ref{Key Lemma} we have 
        \begin{equation}
        \label{Second part}
            \limsup_{t\to 0}\frac{h_t\cos(\dist_\Y(u(\cdot),o))(x)-\cos(\dist_\Y(u(x),o))}{t}\leq -\cos(\dist_{\Y}(u(x),o))|\de u|^2_{\rm HS}.
        \end{equation}
        Combining \eqref{First part} with \eqref{Second part} we finally get the thesis.
\end{proof}
\subsection{A variant of the Bochner-Eells-Sampson inequality}
The authors in \cite{ZZ18} are able to prove the Lipschitz continuity of harmonic maps between Alexandrov spaces exploiting the properties of the Hopf-Lax semigroup. Moreover in \cite{ZZZ19}, given the Lipschitz continuity of the harmonic map proved in \cite{Ser95}, they are able to prove a weak version of the Bochner-Eells-Sampson inequality for maps from a Riemannian domain to a ${\rm CAT}(1)$ space. Here we shall exploit the ideas contained in \cite{Gig23} and fuel them with the ideas of \cite{ZZZ19} (see also \cite{MS22} for the non-smooth counterpart, as in our case) to obtain a variational inequality (the "good" distributional bound) which in the limit will be the desired inequality.

\begin {comment}
\begin{lemma}\label{lem4.1}
 For any $t\in (0,t_*)$, we have the following properties:
 \begin{enumerate}
\item   $f_t\in W^{1,2}_{\rm loc}(B_R)\cap C(B_R)$;
\item  $L_t$ and $D_t$ are lower semi-continuous on $B_R$;
\item  there exist constants $C_1, C_2>0$ depending only $p$ and $A$ such that 
\end{enumerate}
\begin{align*}
    \frac{L_t(x)}{t}&\leq C_{2}\cdot \big[{\rm lip}u(x,C_{1}t^{\frac{p-1}{p}})\big]^{\frac{1}{p-1}};\\
  \frac{D_t(x)}{t}&\leq C_{2}\cdot \big[{\rm lip}u(x,C_{1}t^{\frac{p-1}{p}})\big]^{\frac{p}{p-1}};\\
  \frac {-f_t(x)}{t}&\leq C_{2}\cdot \big[{\rm lip}u(x,C_{1}t^{\frac{p-1}{p}})\big]^{\frac{p}{p-1}}.
\end{align*}
\end{lemma}

\begin{proof}
The first two items were proved in \cite{ZZZ19}.

For item (3), let $y_t\in S_t(x)$ be such that $L_t(x)=\dist_\X(x,y_t)$. Firstly, from $f_t\leq0$ and $F(x,y)\leq A(1+A)$, we get that 
$$\frac{L_t^p(x)}{pt^{p-1}}=F(x,y_t)+f_t(x)\leq A(A+1) \quad \Longrightarrow\quad L_t(x)\leq C_{1}t^{\frac{p-1}{p}},$$
where $C_1:=[pA(A+1)]^{1/p}.$
By the definition of ${\rm lip}u(x,r)$, we have that 
$$D_t(x)=F(x,y_t)\leq (1+A){\dist}_{\Y}\big(u(x),u(y_t)\big)\leq (1+A) \cdot{\rm lip}u(x,C_1t^{\frac{p-1}{p}})\cdot L_t(x).$$
Substituting   this into  
$$\frac{L_t^p(x)}{pt^{p-1}}=D_t(x)+f_t(x)\leq D_t(x),$$
we conclude that 
\begin{align*}
    \frac{L_t (x)}{ t }&\leq \Big[p (1+A) \cdot{\rm lip}u(x,C_1t^{\frac{p-1}{p}})\Big]^{\frac{1}{p-1}};\\
\frac{D_t (x)}{ t }&\leq  (1+A) \cdot{\rm lip}u(x,C_1t^{\frac{p-1}{p}})\cdot \frac{L_t(x)}{t}\leq C_{p,A}\cdot \Big[{\rm lip}u(x,C_1t^{\frac{p-1}{p}})\Big]^{\frac{p}{p-1}}.
\end{align*}

At last, by using $-f_t(x)=D_t(x)-\frac{L^p_t(x)}{pt^{p-1}}\leq D_t(x),$ we obtain the desired upper bound of $-f_t(x)/t$. The proof is now complete.
\end{proof}
\end{comment}

We now recall \cite[Lemma 6.13]{Gig23}.
\begin{lemma}
\label{Perturbation lemma}
    There exists $T>0$ such that, given a Borel set $E\subset B^\prime$ such that $\m(B^\prime\setminus E)=0$, we have: for all $0<t<T$ there exists $z_t\in B$ such that for $\m$-a.e. $x\in E\cap B_{4r/3}(\Bar{x})=:E\cap B^{\prime\prime}$ and every $n\in\N$ the function
    \begin{equation*}
       y\mapsto g_t(x,y,z_t):=\frac{\dist_{\X}^p(x,y)}{pt^{t-1}}+f(x,y)+\frac{\dist^2_\X(y,z_t)}{2n}
    \end{equation*}
    admits a minimizer $T_t(x)$ and such minimizer belongs to the set $E\cap B^{\prime\prime}$.
\end{lemma}
\begin{proof}
    The difference with respect to \cite[Lemma 6.13]{Gig23} lies in the different definition of $f$, however since the proof follows with minor modifications we decided to omit it. Note moreover that from the proof in \cite{Gig23} we can infer that Sobolev regularity is not necessary for the function $f$. It would be sufficient to ask for $f$ to be continuous and with a Laplacian bound $\Delta f\leq L\m$ in the weak sense.
\end{proof}
We further define
\begin{equation}
    \label{f t n}
    f_{t,n}(x):=\inf_{y\in\X}\bigg[\frac{\dist_{\X}^p(x,y)}{pt^{p-1}}+f(x,y)+\frac{\dist^2_\X(y,z_t)}{2n}\bigg].
\end{equation}
We now have the following distributional bound for the function $f_{t,n}$.
\begin{lemma}["Bad" distributional bound]
    Possibly choosing a smaller $t_*$ the following holds. Let $f_{t,n}$ be defined as in \eqref{f t n} and $p\geqslant 2$: we have 
    \begin{equation}
        \label{bad distributional bound}
        {\Delta f_{t,n}}\leq C(K,N,p,{\rm diam}(\Omega))\bigg(\frac{1}{t}+\frac{1}{n}\bigg)\m\qquad{\rm on}\;B
    \end{equation}
    in the weak sense, for all $t<t_*$, for all $n\in\N$.
\end{lemma}
\begin{proof}
Fix $y\in B$: by Theorem \ref{Laplacian comparison} and $p\geqslant 2$, we have 
$$\Delta {\dist}_{\X,y}^p= p(p-1){\dist}_{\X,y}^{p-2}|\nabla \dist_{\X,y}|^2\cdot\m+p\dist_{\X,y}^{p-1}\Delta \dist_{\X,y}\leqslant C(K,N,p,{\rm diam}(\Omega))\cdot\m.$$
Combining this with
 Lemma \ref{Distributional laplacian bound for f} and \cite[Lemma 4.7]{Gig23} we infer the result.
\end{proof}
To obtain the "good" distributional bound we need the following lemma for the function $F$ to be able to let the heat flow and the Hopf-Lax semigroup combine in an efficient way.
\begin{lemma}[Key technical Lemma]
    Consider $4$ points $P,Q,R,S$ inside $u(B_r(x))$ in such a way that $P:=u(x)$, $Q:=u(\Bar{x})$, $R:=u(\Bar{y})$, $S:=u(y)$. Let us further set  $l_0:=2\sin\frac{\dist_\Y(Q,R)}{2}$, $l_1:=2\cos\frac{\dist_\Y(Q,R)}{2}$, $\alpha:=1/(1+2l_0)$ and finally let $\beta>0$. We have
    \begin{equation}
      \label{gennaro}
      F(u(\Bar{x}),u(\Bar{y}))-F(u(x),u(y))\leq\frac{ \bigl[w_{a_1,b,Q_m,\Bar{x}}(x)-w_{a_1,b,Q_m,\Bar{x}}(\Bar{x})\bigr]+\bigl[w_{a_2,b,Q_m,\Bar{y}}(y)-w_{a_2,b,Q_m,\Bar{y}}(\Bar{y})\bigr]}{\alpha l_0},
  \end{equation}
   where $Q_m$ is the middle point of the geodesic joining $Q$ and $R$,  
  \begin{equation*}
      a_1:=1-\frac{1-\alpha}{2}\biggl(1-\frac{1}{\beta}\biggr),\quad b:=l_1,\quad a_2:=1-\frac{1-\alpha}{2}\biggl(1-\beta\biggr)
  \end{equation*}
  and the function $w$ is defined in Lemma \ref{Key Lemma special function}.
\end{lemma}
\begin{proof}
    We can  apply \eqref{technical result} to get  
    \begin{align*}
      \alpha l_0\bigl(F(Q,R)-F(P,S))\bigr)&=\alpha l_0\biggl(4\sin^2\frac{\dist_{QR}}{2}-4\sin^2\frac{\dist_{PS}}{2}\biggr)+\alpha l_0\biggl(2\sin\frac{\dist_{QR}}{2}-2\sin\frac{\dist_{PS}}{2}\biggr) \\
      &\leq\biggl[1-\frac{1-\alpha}{2}\biggl(1-\frac{1}{\beta}\biggr)\biggr]4\sin^2\frac{\dist_{PQ}}{2}+l_1\biggl(\cos\dist_{PQ_m}-\cos\dist_{QQ_m}\biggr) \\ 
      &+\biggl[1-\frac{1-\alpha}{2}\biggl(1-\beta\biggr)\biggr]4\sin^2\frac{\dist_{RS}}{2}+l_1\biggl(\cos\dist_{SQ_m}-\cos\dist_{RQ_m}\biggr) \\
      &\leq \bigl[w_{a_1,b,Q_m,\Bar{x}}(x)-w_{a_1,b,Q_m,\Bar{x}}(\Bar{x})\bigr]+\bigl[w_{a_2,b,Q_m,\Bar{y}}(y)-w_{a_2,b,Q_m,\Bar{y}}(\Bar{y})\bigr],
  \end{align*}
 
  which concludes the proof.
\end{proof}
The second tool we need is an improvement of the earlier distributional bound: this is the aim of the following proposition.
\begin{prop}["Good" distributional bound]
    We have
    \begin{equation}
        \label{good distributional bound}
        {\Delta f_t}_{}\leq -K\frac{L_t^p}{t^{p-1}}+\bigg(1+o_t(1)\bigg)D_t|\de u|^2_{\rm HS}\quad{\rm on}\; B
    \end{equation}
    in the weak sense, for all $t<t^*$ and $p\geqslant2$.
\end{prop}
\begin{proof}
    First of all let us recall the definition of $f_{t,n}$ 
    \begin{equation*}
    f_{t,n}(x):=\inf_{y\in\X}\bigg[\frac{\dist^p_\X(x,y)}{pt^{p-1}}+f(x,y)+\frac{\dist^2_\X(y,z_t)}{2n}\bigg].
    \end{equation*}
    Thanks to the Lemma \ref{Perturbation lemma} we can find $z_t$ in such a way that a minimizer of $g_t(x,y,z_t)$, i.e. a point $T_t(x)$ for which $g_t(x,T_t(x),z_t)=f_{t,n}(x)$, lies inside $E\cap B^{\prime\prime}$ for $\m$-a.e. $x\in E\cap B^{\prime\prime}$ and we can choose $E$ to be the set of regular points of the space intersected with the set of Lebesgue points of $|\de u|_{\rm HS}$ (which is clearly of full measure). Now let us fix $\Bar{x}\in E\cap B^{\prime\prime}$ and call $\Bar{y}$ the "good" minimiser of $f_{t,n}(x)$. Clearly for such points we have 
    \begin{equation*}
        f_{t,n}(\Bar{x})=f(\Bar{x},\Bar{y})+\frac{\dist_{\X}^p(\Bar{x},\Bar{y})}{pt^{p-1}}+\frac{\dist_{\X}^2(\Bar{y},z_t)}{2n}=-F(u(\Bar{x}),u(\Bar{y}))+\frac{\dist_{\X}^p(\Bar{x},\Bar{y})}{pt^{p-1}}+\frac{\dist_{\X}^2(\Bar{y},z_t)}{2n}.
    \end{equation*}
    Now fix any other two points $x,y\in\Omega$. Setting  $P:=u(x)$, $Q:=u(\Bar{x})$, $R:=u(\Bar{y})$, $S:=u(y)$. Using the inequality \eqref{gennaro} of the key techcnical lemma (and its notation) we get
\begin{align*}
    f_{t,n}(x)&=\inf_{y\in\X}\bigg[\frac{\dist^p_\X(x,y)}{pt^{p-1}}+f(x,y)+\frac{\dist^2_\X(y,z_t)}{2n}\bigg] \\
    &=-F(u(\Bar{x}),u(\Bar{y}))+\inf_{y\in\X}\bigg[\frac{\dist^p_\X(x,y)}{pt^{p-1}}+F(u(\Bar{x}),u(\Bar{y}))-F(u(x),u(y))+\frac{\dist^2_\X(y,z_t)}{2n}\bigg] \\
    &\underbrace{\leq}_{\eqref{gennaro}}-F(u(\Bar{x}),u(\Bar{y}))+\frac{w_{a_1,b,Q_m,\Bar{x}}(x)-w_{a_1,b,Q_m,\Bar{x}}(\Bar{x})}{\alpha l_0}\\
    &\qquad+Q_t\bigg[\frac{w_{a_2,b,Q_m,\Bar{y}}(\cdot)-w_{a_2,b,Q_m,\Bar{y}}(\Bar{y})}{\alpha l_0}+\frac{\dist_\X^2(\cdot,z_t)}{2n}\bigg](x),
\end{align*}
with equality if $x=\Bar{x}$.
We now proceed to obtain a bound on the Laplacian of $f_{t,n}$ in the heat flow sense at the point $\Bar{x}$, therefore we shall estimate
\begin{equation*}
    \limsup_{s\to 0^+}\frac{h_s(f_{t,n})(\Bar{x})-f_{t,n}(\Bar{x})}{s}=\Delta f_{t,n}(\Bar{x}).
\end{equation*}
Exploiting the previous inequalities and the monotonicity of the heat flow ($h_tf\leq h_t g$ if $f\leq g$) we get 
\begin{equation*}
    \Delta f_{t,n}(\Bar{x})\leq \frac{\Delta w_{a_1,b,Q_m,\Bar{x}}(\Bar{x})}{\alpha l_0}+\Delta Q_t\bigg[\frac{w_{a_2,b,Q_m,\Bar{y}}(\cdot)}{\alpha l_0}+\frac{\dist_\X^2(\cdot,z_t)}{2n}\bigg](\Bar{x}).
\end{equation*}
Moreover thanks to the properties of the Hopf-Lax semigroup (namely \eqref{Heat flow bound Hopf-Lax}) we get 
\begin{equation*}
    \Delta Q_t\bigg[\frac{w_{a_2,b,Q_m,\Bar{y}}(\cdot)}{\alpha l_0}+\frac{\dist_\X^2(\cdot,z_t)}{2n}\bigg](\Bar{x})\leq \frac{\Delta w_{a_2,b,Q_m,\Bar{y}}(\Bar{y})}{\alpha l_0}+\frac{1}{n}\Delta\dist_\X^2(\cdot,z_t)(\Bar{y})-K\frac{L_t^p(\Bar{x})}{t^{p-1}}.
\end{equation*}
Now we can apply Lemma \ref{Key Lemma special function} and the Laplacian comparison to obtain 
\begin{align*}
    \Delta f_{t,n}(\Bar{x})&\leq \frac{C(K,N,r)}{n}-K\frac{L_t^p(\Bar{x})}{t^{p-1}}+\frac{2a_1-b\cos(\dist_\Y(u(\Bar{x}),Q_m))}{\alpha l_0}|\de u|^2_{\rm HS}(\Bar{x}) \\ 
    &+\frac{2a_2-b\cos(\dist_\Y(u(\Bar{y}),Q_m))}{\alpha l_0}|\de u|^2_{\rm HS}(\Bar{y}).
\end{align*}
 Since $\cos(\dist_\Y(Q_m),\Bar{y})=\cos(\dist_\Y(Q_m),\Bar{x})=l_1/2$ and $1-l_1^2/4=l_0^2/4$ we can choose $\beta$ such that 
  $a_2=l_1^2/4$, so that $2a_2-b\cos(\dist_\Y(u(\Bar{y}),Q_m))=0$. This is achieved with  
  \begin{equation*}
      \beta = 1-\frac{l_0(1+2l_0)}{4}.
  \end{equation*}
  Via standard computations we get 
  \begin{equation*}
      \frac{2a_1-b\cos(\dist_\Y(Q_m),\Bar{x})}{\alpha l_0}=2l_0(1+2l_0)\bigg(\frac{1}{4}+\frac{1}{4-l_0(1+2l_0)}\bigg).
  \end{equation*}
  Therefore we get 

\begin{align*}
     \Delta f_{t,n}(\Bar{x})&\leq \frac{C(K,N,r)}{n}-K\frac{L_t^p(\Bar{x})}{t^{p-1}}+2l_0(1+2l_0)\bigg(\frac{1}{4}+\frac{1}{4-l_0(1+2l_0)}\bigg)|\de u|^2_{\rm HS}(\Bar{x}) \\ 
     &\leq \frac{C(K,N,r)}{n}-K\frac{L_t^p(\Bar{x})}{t^{p-1}}+\bigg(1+o_t(1)\bigg)D_t(\Bar{x})|\de u|^2_{\rm HS}(\Bar{x}).
\end{align*}
where we also used that $D_t(\Bar{x})=l_0+l_0^2$ and that $u$ is H\"older continuous to estimate the remainder in $o_t(1)$  (observe also that $\hat{x}$ does not depend on $n\in\N$).
Combining the latter with Lemma \ref{bad distributional bound} and \cite[Lemma 4.8]{Gig23} (recalling that $u$ is continuous on $\Omega$) we end up with
\begin{equation*}
      \Delta f_{t,n}\leq  \frac{C(K,N,r)}{n}-K\frac{L_t^p(\cdot)}{t^{p-1}}+\bigg(1+o_t(1)\bigg)D_t(\cdot)|\de u|^2_{\rm HS}(\cdot)\quad{\rm on}\; B
\end{equation*}
in the weak sense, for all $n\in\N$ and for all $t<t_*$ .

Now since $f_{t,n}$ converges to $f_t$ uniformly as $n\to\infty$, thanks to the regularity of $f_t$ and the stability of the Laplacian bounds we infer \eqref{good distributional bound}.

\end{proof}

We now recall \cite[Lemma 4.4]{ZZZ19}:
\begin{lemma}
\label{Limits for Bochner}
    Let $q$ be such that $1/q+1/p=1$. For all $x\in B$ we have 
    \begin{equation}
    \label{Limit Bochner 1}
        \liminf_{t\to 0}\frac{f_t(x)}{t}\geq-\frac{1}{q}\lip^qu(x).
    \end{equation}
    Moreover, in adding   to assume that $u$ is locally Lipscithz continuous,  for $\m$-a.e. $x\in B$ (namely any point in $B$ where $u$ is metrically differentiable) we have
    \begin{equation}
        \label{Limit Bochner 2}
        \lim_{t\to 0^+}\frac{f_t(x)}{t}=-\frac{{\rm lip}^qu(x)}{q}
    \end{equation}
    and
    \begin{equation}
        \label{Limit Bochner 3}
        \lim_{t\to 0^+}\frac{L_t(x)}{t}={\rm lip}^{q/p}u(x),\qquad\lim_{t\to 0^+}\frac{D_t(x)}{t}={\rm lip}^qu(x).
    \end{equation}
\end{lemma}
\begin{proof}
    The proof follows as in \cite[Proposition 7.5]{MS22} combined with   \cite[Lemma 4.4]{ZZZ19}.
\end{proof}

\begin{teo}[A variant of the BES inequality]\label{thm:varbes}
Let $u$ be as above and assume that it is locally Lipschitz in $\Omega$, then the inequality  
\begin{equation}
    \label{Bochner inequality}
    \Delta\bigg(\frac{\lip^2u}{2}\bigg)\geq|\nabla\lip u|^2-K\lip^2(u)- e_2^2[u]\lip^2u
\end{equation}
holds in the very weak sense in $\Omega$.
\end{teo}
\begin{proof}
     By the chain rule it is easy to infer that \eqref{Bochner inequality} is equivalent to 
    \begin{equation}
        \label{Reduced Bochner}
        \Delta\lip u\geq-K\lip u- e_2^2[u]\lip u.
    \end{equation}
We shall now verify that there exists a neighborhood $B_R(\Bar{x})$ with $B_{2R}(\Bar{x})\subset\Omega$ such that $\lip(u)\in W^{1,2}(B_R(\Bar{x}))$ and \eqref{Reduced Bochner} holds in the sense of distributions in $B_R(\Bar{x})$.

Due to the continuity of $u$ there exists $R>0$ such that $u(B_{2R}(\Bar{x}))\subset B_{\pi/4}(u(\Bar{x}))$, so that ${\rm diam}(u(B_{2R}(\Bar{x})))<\pi/2$ and $R<r/2$. By \eqref{good distributional bound} and \eqref{Modified Lemma} we have $\Delta f_t/t\leq C(\Lip u)$ on $B_{2R}$ for all $t\in(0,t_*)$. Combining the elliptic inequality \eqref{good distributional bound} with Lemma \ref{Modified Lemma} and a Caccioppoli inequality we get $f_t/t\in W^{1,2}(B_{3R/2}(x))$ with $\|f_t/t\|_{W^{1,2}(B_{3R/2}(\Bar{x}))}\leq C$ and $C$ depending only on the Lipschitz norm of $u$ in $\Omega^\prime$. Therefore, exploiting Lemma \ref{Limits for Bochner}, up to a subsequence we have that $-f_t/t$ converges weakly in $W^{1,2}$ to $\lip^q(u)/q$ and we get
\begin{equation}
\label{Almost Bochner}
    \Delta(\lip^qu/q)\geq K\lip^q u- e_2^2[u]\cdot\lip^q u
\end{equation}
in $B_{3R/2}(\Bar{x})$ in the weak sense.
Exploiting the Lipschitz continuity of $u$  we get 
\begin{equation*}
    \Delta(\lip^qu/q)\geq K(\lip u)^q-(\lip u)^{q+2}\geq -C
\end{equation*}
where the constant is uniform in $q$. Now again by Caccioppoli inequality we get $\|\lip^q u/q\|_{W^{1,2}(B_R)}\leq C$ as $q\to 1$. This means that $\lip^q(u)/q$ converges to $\lip(u)$ in $W^{1,2}(B_R(\Bar{x}))$ and we can pass to the limit in \eqref{Almost Bochner} and get \eqref{Reduced Bochner}, whence we also deduce \eqref{Bochner inequality}. 
\end{proof}
Finally we shall mention that the theorems in \cite[Section 5]{ZZZ19} hold also in the present setting: we refer to \cite{ZZZ19} for the proofs which work mutatis mutandis in our context.
\begin{teo}\label{thm:sharplip}
    Let $u$ be as above but with values in $B_\rho(o)\subset\Y$, where $(\Y,\dist_\Y)$ is a ${\rm CAT}(\kappa)$ space and $\rho<\pi/2\sqrt{\kappa}$. Then letting $R>0$ be such that $B_{2R}(x_0)\subset\Omega$ we have
    \begin{equation}
    \label{Yau inequality}
        \sup_{x\in B_{R/2}(x_0)}{\rm lip}(u)(x)\leq\frac{C_{N,\sqrt{K}R,\pi/(2\sqrt{\kappa}-\rho)}}{R},
    \end{equation}
    where the constant $C$ only depends on the parameters listed at its subscript.
\end{teo}
As a consequence we obtain a Liouville type theorem for harmonic maps, which follows by estimate \eqref{Yau inequality}.
\begin{cor}\label{Yauineq}
    Let $(\X,\dist_\X,\m_X)$ be an ${\rm RCD}(0,N)$ space and $(\Y,\dist_\Y)$ be a ${\rm CAT}(\kappa)$ space. Consider an harmonic map $u:\X\to\Y$ such that $u(\X)\subset B_\rho(o)$ for some $o\in\Y$ and $\rho<\pi/(2\sqrt{\kappa})$ with sublinear growth, i.e.
    \begin{equation*}
        \liminf_{R\to\infty}\frac{\sup_{y\in B_R(x_0)}\dist_\Y(u(y),o)}{R}=0
    \end{equation*}
    for some $o\in\Y$. Suppose that $u$ is locally Lipschtiz continuous. Then $u$ must be a constant map.
\end{cor}

 \subsection{Boundary regularity for harmonic maps}
 
In this section, we continue to assume that  $\Omega\subset X$ is an open bounded set in an ${\rm RCD}(K,N)$ space with $X\setminus \Omega\not=\emptyset$ , $K\in\mathbb R$ and $N\in[1,\infty)$.  Moreover we let $(Y,{\dist}_{\Y})$ be a ${\rm CAT}(\kappa)$ space with $\kappa>0$. 

To study the boundary regularity of harmonic maps, we shall also impose some regularity conditions on the boundary of $\Omega.$
\begin{defn}
    Let $\Omega\subset X$ be a  domain. We say that $\Omega $ satisfies an {\emph{exterior density condition}} if there exist two numbers $\lambda\in(0,1)$ and $R_0>0$ such that  
\begin{equation} \label{equ-1.3}
\m(\Omega\setminus B_r(x))\geqslant \lambda \cdot \m(B_r(x)) \quad \forall\ x\in \partial\Omega,\ \  \forall r\in(0,R_0).
\end{equation}
Additionally we say that $\Omega$ satisfies a {\emph {uniform exterior sphere condition}} if there exists a number $R_0>0$ such that for each $x_0\in\partial\Omega$ there exists a ball
 $ B_{R_0}(y_0)$ satisfying
\begin{equation}\label{equ-1.4}
\Omega\cap B_{R_0}(y_0)=\emptyset\quad {\rm and}\quad x_0\in \partial B_{R_0}(y_0).
\end{equation}
\end{defn}
\begin{remark}
    It is easy to see that if the space satisfies a volume doubling condition (which is the case of ${\rm RCD}(K,N)$ spaces, thanks to Bishop-Gromov inequality), then the exterior density condition is implies by the exterior sphere condition.
\end{remark}
 
The main result of this section is the following.

 \begin{teo} \label{thm-boundary-reg}
Let $\Omega $ and $Y$ be as above. Suppose that $\Omega\subset X$ satisfies a uniform exterior sphere condition with constant $R_0$ and let $w\in {\rm Lip}(\overline{\Omega},Y)$. Let $u\in {\rm KS}^{1,2}(\Omega,Y)$ be an harmonic map with boundary data $w$ such that ${\rm Im}(u)\subset B_{\pi/4-\rho}(o) $ for some $o\in Y$ and $\rho>0$.  Then for any $\epsilon\in(0,1)$ it holds 
\begin{equation}\label{equ-4.15}
 {\dist}_{\Y}\big(u(x),w(x_0)\big)\leqslant C_{\epsilon}   L_w   {\dist}_{\X}^{1-\epsilon}(x,x_0) 
 \end{equation}
 for all  $x_0\in\partial\Omega$ and $x\in\Omega$ with ${\dist}_{X}(x,x_0)<R_\epsilon$, where both $R_\epsilon$ and $C_\epsilon$ depend only on $\epsilon, N,K$ and ${\rm diam}(\Omega)$, and  
  $$L_w:= \sup_{x,y\in\overline\Omega} \frac{{\dist}_{\Y}\big(w(x),w(y)\big)}{{\dist}_{X}(x,y)}.$$ 
  In particular, $u$ is continuous at $x_0$ and $u(x_0)=w(x_0)$.
 \end{teo}

To prove this result, we need the following two lemmas.

 % \begin{teo}\label{boundary-reg}
 %Let $\Omega $ and $ Y$ and $u$ be as in the above. Let $w\in {\rm KS}^{1,2}(\Omega,Y)$ . %Suppose that $u\in{\rm KS}(\Omega;Y)$ be a harmonic map with  is the boundary data of $u$, i.e., ${\dist}_{\Y}(u(x),w(x))\in  W^{1,2}_0(\Omega)$  and 
 %with values in a ball $B_\rho(o)\subset Y$ for some $\rho\in(0,\frac{\pi}{2\sqrt{\kappa}})$.   %We shall always fix a H\"older continuous representative of $u$.
 %Then:
 
%(A) If  $\partial\Omega$ satisfies an exterior density condition and   if $w\in C^\alpha({\overline\Omega},Y)$ and $u\in C^\alpha_{\rm loc}(\Omega,Y)$, then $u\in C^\beta(\overline{\Omega},Y)$ for some small $  \beta\in(0,\alpha)$.
 
% (B) If $\partial \Omega$ satisfies a uniformly exterior sphere condition  if $w\in Lip({\overline\Omega},Y)$, then $u\in C^{1-\epsilon}(\overline{\Omega},Y)$ for every $\epsilon>0$. 
 % \end{teo} 

\begin{lemma}\label{boundary-reg-harmonic-function}
Let $\Omega\subset X$ be a bounded domain satisfying a uniformly exterior condition with constant $R_0$. Suppose that $f\in   W^{1,2}(\Omega)$ is a harmonic function on $\Omega$ with boundary data $g\in {\rm Lip}(\overline{\Omega})$. Suppose $g(z_0)=0$ for some $z_0\in\overline\Omega$. Then for any $\epsilon\in(0,1)$, there exists a number $R_\epsilon\in (0,\min\{1,R_0/2\})$ (depending only on $\epsilon, N, K$ and ${\rm diam}(\Omega)$) such that for any ball $B_r(x_0)$ with  $x_0\in\partial \Omega$  and $r\in (0,R_\epsilon)$ it holds 
\begin{equation}\label{equ-4.12}
 \sup_{B_r(x_0)\cap \Omega} |f(x)-f(x_0)|\leqslant C_{\epsilon}L\cdot r^{1-\epsilon}, 
 \end{equation}
  where the constant $C_{\epsilon}>0$ depending only on $\epsilon, N,K$, and the constant $L$ is a Lipschitz constant of $g$.  

\end{lemma}
 \begin{proof}
This  is Theorem 4.3 in \cite{zz2024}.   
  \end{proof}

\begin{lemma}\label{lem-sub-harmonic} 
Let $\Omega, Y $ be as above.  Suppose that $u:\Omega\to Y$ is an harmonic map. Then  for any 
$P\in Y$ such that ${\rm Im}(u)\in B_{\pi/2-\rho}(P)$  it holds
  \begin{equation}\label{equ-sub-har}
  { \Delta} {\dist}_{\Y}\big(u(x),P\big)\geqslant 0 
  \end{equation}
  in the sense of distributions.  
  \end{lemma}
\begin{proof}
Since the function ${\dist}_{\Y}(P,\cdot)$ is convex in $B_{\pi/2}(P)\subset Y$, the assertion follows directly from  Proposition \ref{convex + harmonic}.
\end{proof}

We are now in the position to prove Theorem \ref{thm-boundary-reg}, whose proof is a modification of the one in \cite[Theorem 4.6]{zz2024}.

 \begin{proof}[Proof of Theorem \ref{thm-boundary-reg}]
 Fix any a point $x_0\in \partial\Omega$, and set  $P=w(x_0).$ Then, by the triangle inequality and the fact that ${\rm Im}(u)\subset B_{\pi/4-\rho}(o)$, we have $ {\dist}_{\Y}(P,u(x))\leqslant \pi/2-2\rho$ for any $x\in \Omega$. Moreover by Lemma \ref{lem-sub-harmonic}, we observe that ${\dist}_{\Y}(P,u(x))$ is sub-harmonic on $\Omega$. 
 
 We can now solve the Dirichlet problem 
 $${\Delta} f(x) =0\quad {\rm on}\ \Omega\qquad {\rm and}\qquad f(x)-{\dist}_{\Y}(w(x_0),w(x))\in W^{1,2}_0(\Omega).$$
 Notice that, by the triangle inequality, the function $g_{x_0}(x) :={\dist}_{\Y}(w(x_0),w(x))$ is Lipschitz continuous on $\overline\Omega$ with a Lipschitz constant 
 $$L_{g_{x_0}}  \leqslant  L_w\quad {\rm and}\quad g_{x_0}(x_0)=0.$$
 According to Lemma  \ref{boundary-reg-harmonic-function}, we have 
\begin{equation} \label{equ-boundary-har-func} 
 \sup_{B_r(x_0)\cap\Omega} |f(x)-f(x_0)|\leqslant C_{\epsilon}  L_w  r^{1-\epsilon}, 
 \end{equation}
   for any ball $B_r(x_0)$ with  $x_0\in\partial \Omega$  and $r\in (0,R'_\epsilon)$. 
   
 At last, since ${\dist}_{\Y}(u(x),w(x_0))-f(x)$ is sub-harmonic on $\Omega$, and     $$[{\dist}_{\Y}(u(x),w(x_0))-f(x)]^+\in W^{1,2}_0(\Omega),$$
the maximum principle yields
   $${\dist}_{\Y}(u(x),w(x_0))\leqslant f(x),\quad {\rm a.e.\ in }\ \ \Omega.$$ 
  Noticing that $u\in C(\Omega)$ (by Theorem \ref{Holder continuity}) and $f\in C(\Omega)$, we get  
   $${\dist}_{\Y}(u(x),w(x_0))\leqslant f(x),\quad  \forall x\in  \Omega.$$  
The combination of the latter with (\ref{equ-boundary-har-func}) implies the desired result, concluding the proof.
   \end{proof}

\bibliography{references}
\bibliographystyle{alpha}
\end{document}